\documentclass{m2an}
\usepackage{graphicx}
\newcommand{\re}{{\rm Re}\,}

\newtheorem{ssmptn}[thrm]{Assumption}
\begin{document}
%%%%%%%%%%%%%%%%%%%%%%%%%%%%%%%%%%%%%%%%%%%%%%%%%%%%%%%%%%%%%%%%%%%%%%%%
% from old sec_2-1.tex
%%%%%%%%%%%%%%%%%%%%%%%%%%%%%%%%%%%%%%%%%%%%%%%%%%%%%%%%%%%%%%%%%%%%%%%%
\title{Initial-boundary value problems for second order systems of partial
differential equations}\thanks{This work performed under the auspices of the
U.S. Department of Energy by Lawrence Livermore National Laboratory under
Contract DE-AC52-07NA27344. O.E.O.~acknowledges support by grants 05/B415 and
214/10 from SeCyT-Universidad Nacional de C\'ordoba, 11220080100754 from
CONICET, PICT17-25971 from ANPCYT, and the Partner Group grant of the Max Planck
Institute for Gravitational Physics, Albert-Einstein-Institute (Germany).}
\author{Heinz-Otto Kreiss}\address{Tr\"ask\"o-Stor\"o Institute of Mathematics, Stockholm, Sweden.}
\author{Omar E. Ortiz}\address{Facultad de Matem\'atica, Astronom\'\i a y F\'\i
sica, Universidad Nacional de C\'ordoba and IFEG, Argentina.}
\author{N. Anders Petersson}\address{Center for Applied Scientific
  Computing, Lawrence Livermore National Laboratory, Livermore, California, USA. }
\date{August 27th, 2009, revised August 12th, 2010.}

\begin{abstract}
We develop a well-posedness theory for second order systems in bounded domains where
boundary phenomena like glancing and surface waves play an important role.  Attempts have
previously been made to write a second order system consisting of $n$ equations as a
larger first order system. Unfortunately, the resulting first order system consists, in
general, of more than $2n$ equations which leads to many complications, such as 
side conditions which must be satisfied by the solution of the larger first order
system. Here we will use the theory of pseudo-differential operators combined with mode
analysis. There are many desirable properties of this approach: 1) The reduction to first
order systems of pseudo-differential equations poses no difficulty and always gives a
system of $2n$ equations. 2) We can localize the problem, i.e., it is only necessary to
study the Cauchy problem and halfplane problems with constant coefficients. 3) The class
of problems we can treat is much larger than previous approaches based on ``integration by
parts''. 4) The relation between boundary conditions and boundary phenomena becomes
transparent.
\end{abstract}
%\begin{resume}
%[ADD FRENCH TRANSLATION OF ABSTRACT]
%\end{resume}

\maketitle

\section*{Introduction}

The theory for first order hyperbolic systems, which was developed with fluid problems in
mind, is by now rather well understood. It turned out that energy estimates via
`integration by parts' and characteristics are the most important ingrediencies in the
theory.

Second order hyperbolic systems often describe problems where wave propagation is
dominant. In bounded domains this leads to a large number of boundary phenomena like
glancing waves and surface waves. Attempts have previously been made to write a second
order system consisting of $n$ equations as a larger first order system. However, boundary
phenomena such as glancing and surface waves correspond to generalized eigenvalues which
are not handled by the theory for first order systems. Furthermore, the resulting first
order system often consists of more than $2n$ equations which leads to many
complications. In particular, the first order system must in general be augmented by side
conditions to guarentee that solutions of the first order system satisfy the original
second order system.

In this paper we describe a theory for second order hyperbolic systems based on Laplace
and Fourier transform, with particular emphasis on boundary processes corresponding to
generalized eigenvalues. Our theory uses pseudo-differential operators combined with mode
analysis, and builds upon the theory for first order systems developed
in~\cite{Kreiss1970,Kreiss-Lorenz1989} . This approach has many desirable properties: 1)
Once a second order system has been Laplace and Fourier transformed it can always be
written as a system of $2n$ first order pseudo-differential equations. Therefore, the
theory of \cite{Kreiss1970,Kreiss-Lorenz1989} also applies here. 2) We can localize the
problem, i.e., it is only necessary to study the Cauchy problem and halfplane problems
with constant coefficients. 3) The class of problems we can treat is much larger than
previous approaches based on ``integration by parts''. 4) The relation between boundary
conditions and boundary phenomena becomes transparent.

The remainder of the paper is organized as follows. In section \ref{sec_2} we state the general
problem and provide some basic definitions. In section \ref{sec_3} we treat in detail the
fundamental problem of a single wave equation in a half-plane subject to different types of boundary
conditions.  In section \ref{sec_4} we first study two wave equations coupled through the boundary
conditions and then outline a theory for the general case of $n$ second order wave equations. This
theory proves that all essential difficulties already occur for scalar wave equations coupled
through the boundary conditions.  Numerical experiments are presented in section~\ref{sec_5}, where
we study the different classes of boundary phenomena for two wave equations coupled through the boundary conditions.

\section{Initial-Boundary Value Problems for second order hyperbolic
systems}\label{sec_2}
\subsection{Well posed problems}\label{sec_2-1}
In this paper we want to consider second order systems which are of the form
\begin{equation}\label{2.1.1}
u_{tt}=P_0(D) u + F(x,t), \quad t\ge 0, \quad x\in \Omega, \quad F\in
\xCinfty_0(\Omega),
\end{equation}
in the halfspace $\Omega=\{x_1\ge 0, -\infty < x_j <\infty, j=2,\dots,r\}.$ Here
\begin{equation}\label{2.1.2}
P_0(D)= A_1D_1^2+\sum_{j=2}^r B_jD_j^2, 
\end{equation}
where 
\[ A_1=A_1^*>0, ~B_j=B_j^*>0, \]
are $n\times n$ constant matrices, $u$ is a vector valued function with $n$
components and we are using the notation
\[\begin{split}
x&=(x_1,\ldots,x_r),\quad D=(D_1,\ldots,D_r),~D_j=\partial/\partial x_j,\\
u_t &=\partial u/\partial t =D_t u,\quad u_{x_j}=D_j u.
\end{split}
\]
At $t=0$ we give initial conditions by
\[
u(x,0)=f_1(x), \quad u_t(x,0)=f_2(x).
\]
We are interested in smooth solutions which belong to $\xLtwo(\Omega)$ and satisfy,
at the boundary $\Gamma = \{x_1=0, -\infty < x_j <\infty, j=2,\dots,r\}$ $n$
linearly independent boundary conditions.
\begin{equation}
Lu=:C_0 u_t+\sum_{j=1}^r C_j u_{x_j}=g,\quad g\in
\xCinfty_0(\Gamma).\label{2.1.3}
\end{equation}
Here $C_0,C_j$ are constant $n\times n$ matrices, $C_1$ is non-sigular and,
without loss of generality, we assume it to be normalised 
\begin{ssmptn}\label{assumption2.1.1} $C_1 = I.$ \end{ssmptn}

To facilitate the use of Laplace transformation in time, we frequently assume
that the initial data are homogeneous, i.e., $f_1=f_2\equiv 0.$ This is however
no restriction, since it is always possible to change variables in a problem
with general initial data such that the initial data becomes homogeneous in the
new variable. Since the Cauchy problem is well posed (see section \ref{sec_5-1})
we can extend the definition of the forcing and the initial data smoothly to the
whole of ${\mathbb R}^r(x)$ and determine its solution. Then we subtract this
solution from the halfplane problem and obtain a new halfplane problem where
only the boundary data do not vanish. This is a very natural procedure because
all the difficulties and many physical phenomena arise at the boundary.

We now introduce some key definitions that classify the problems according to estimates one can achieve.

\begin{dfntn}\label{definition2.1.1}
Consider (\ref{2.1.1})--(\ref{2.1.3}) for $F=0,$ $f_1=f_2=0.$
The problem is called {\it Strongly Boundary Stable} if
there are constants $\eta_0\ge 0,$ and  $K>0$ which are
independent of $g,$ such that for all $\eta \ge \eta_0\ge 0,$ $T\ge 0$
\begin{equation}
\int_0^T e^{-2\eta t} \Bigl(\|u(\cdot,t)\|^2_{\xHone(\Gamma)} +
\|u_t(\cdot,t)\|^2_{\xHzero(\Gamma)}\Bigr)~dt \le K 
~\int_0^T e^{-2\eta t} \|g(\cdot,t)\|^2_{\xHzero(\Gamma)}
~dt.\label{2.1.4}
\end{equation}
\end{dfntn}
\begin{dfntn}\label{definition2.1.2}
The problem (\ref{2.1.1})--(\ref{2.1.3}) is called {\it Boundary Stable} if
there are constants $\eta_0>0,$ $K>0$ and $\alpha > 0,$ which are
independent of $g,$ such that for all $\eta \ge \eta_0,~T\ge 0,$
\begin{equation}
\int_0^T e^{-2\eta t} \|u(\cdot,t)\|^2_{\xHzero(\Gamma)}~dt \le \frac{K}{
\eta^\alpha}~\int_0^T e^{-2\eta t} \|g(\cdot,t)\|^2_{\xHzero(\Gamma)}~dt.
\label{2.1.5}
\end{equation}
\end{dfntn}
Here $\|u\|^2_{\xHn{p}}$ denotes the norm composed of the $\xLtwo$-norm of $u$ and all its
derivatives up to order $p.$ Thus (\ref{2.1.4}) tells us that we ``gain'' one derivative
while (\ref{2.1.5}) says that $u$ is as smooth as the data. The constants $\alpha,\eta_0$
are very important. If $\eta_0=0,$ then we can choose $\eta=\frac{1}{T}$ for every fixed
$T>0.$ This shows that the solution grows at most like $T^\alpha$ with time.  If
$\eta_0>0,$ then there is bounded exponential growth. This can happen when lower order
terms are present.

The boundary estimates allow us also to obtain interior estimates. In section
\ref{sec_3-2} we will prove
\begin{thrm}\label{theorem2.1.1}
Consider  (\ref{2.1.1})--(\ref{2.1.3}) with $F=0.$ If the problem is 
{\it Boundary Stable}, then we obtain interior estimates of the form (\ref{2.1.4}),
(\ref{2.1.5}) where $\|u(\cdot,t)\|^2_{\xHzero(\Gamma)}$ is replaced by
$\|u(\cdot,t)\|^2_{\xHzero(\Omega)}$ and $\alpha$ by $\tilde \alpha\ge \alpha +1,$
respectively. If the problem is not {\it Boundary Stable}, then it is illposed.
\end{thrm}

Since we can always reduce the data such that only $g\ne 0,$ we could restrict 
ourselves to this case. However, we are interested in differential equations with
variable coefficients in general domains. Thus we have also to discuss the case that
 $F\ne 0.$ In particular, we have to show that the problem is stable against
perturbations by lower order (first order) terms of the differential equations.

\begin{dfntn}\label{definition2.1.3}
The problem (\ref{2.1.1})--(\ref{2.1.3}) with $f_1=f_2=0$ is called {\it Strongly Stable} if there
exists $\eta_0>0,$ $T>0,$ $K>0$ and $\alpha > 0,$ which are independent of $g$ and $F$ such that,
for all $\eta \ge \eta_0,$
\begin{align}
 \eta \int_0^T & e^{-2\eta t} \Bigl(\|u(\cdot,t)\|^2_{\xHone(\Gamma)} +
\|u_t(\cdot,t)\|^2_{\xHzero(\Gamma)}\Bigr)~dt + \eta^2\int_0^T e^{-2\eta t}
\Bigl(\|u(\cdot,t)\|^2_{\xHone(\Omega)} +
\|u_t(\cdot,t)\|^2_{\xHzero(\Omega)}\Bigr)~dt\nonumber\\
& \le K ~\Bigl[\eta \int_0^T e^{-2\eta t}
\|g(\cdot,t)\|^2_{\xHzero(\Gamma)} ~dt + \int_0^T e^{-2\eta t}
\|F(\cdot,t)\|^2_{\xHzero(\Omega)}~dt\Bigr]. \label{2.1.6}
\end{align}
\end{dfntn}

Clearly, if (\ref{2.1.6}) holds, then the problem is {\it Strongly Boundary Stable}. For first order
systems the classical theory (see \cite{Kreiss1970,Kreiss-Lorenz1989}) tells us that also the
converse is true: If the problem is {\it Strongly Boundary Stable}, then it is {\it Strongly
  Stable}. As we will see, after Laplace and Fourier transformation we can write our problem again
as a first order system which satisfies all the conditions of the classical theory and therefore the
results of that theory are also valid for second order systems.  In particular, the problem is
stable against lower order perturbations both for the differential equations and the boundary
conditions.  (See Appendix of \cite{Kreiss2006}).

Due to physical phenomena like glancing and surface waves, the problems for
second order systems are often only {\it Boundary Stable}. This leads to
\begin{dfntn}\label{definition2.1.4}
We call the problem (\ref{2.1.1})--(\ref{2.1.3}) {\it Stable} if it is {\it Boundary Stable} and if,
 for $g=0,$ there exists
$\eta_0\ge 0,$ $K>0$ and $\alpha > 0$
 which are independent of $F$ such that, for all $\eta \ge \eta_0,$
\begin{equation}\label{2.1.7}
\int_0^T e^{-2\eta t} \Bigl(\|u(\cdot,t)\|^2_{\xHone(\Omega)} +
\|u_t(\cdot,t)\|^2_{\xHzero(\Omega)}+
\|u(\cdot,t)\|^2_{\xHzero(\Gamma)}\Bigr)~dt \le \frac{K}{
\eta^\alpha}~\int_0^T e^{-2\eta t} \|F(\cdot,t)\|^2_{\xHzero(\Omega)}~dt.
\end{equation}
\end{dfntn}
 
If (\ref{2.1.7}) holds, then we can obtain an estimate even when $g\ne 0.$ We split the
problem into two; one with $g=0$ and $F\ne 0$ and the other with $g\ne 0$ and
$F=0.$ For the first problem we obtain (\ref{2.1.7}) and for the other we use
Theorem \Rref{theorem2.1.1}. 

In applications there is often a standard energy estimate, which can be obtained by
integration by parts provided that $g=0$. This estimate can be written as
\begin{equation}\label{eq:energy-est}
\|u(\cdot,t)\|^2_{\xHone(\Omega)} + \|u_t(\cdot,t)\|^2_{\xHzero(\Omega)} \le 
K \left[ \|u(\cdot,0)\|^2_{\xHone(\Omega)} + \|u_t(\cdot,0)\|^2_{\xHzero(\Omega)} +
\int_0^t \|F(\cdot,\tau)\|^2_{\xHzero(\Omega)}\,d\tau \right],
\end{equation}
where the constant $K$ is independent of $F$. In this case we need only to show that the
problem is {\it Boundary Stable}.
\begin{thrm}\label{theorem:1.7}
The problem is {\it Stable} if it is {\it Boundary Stable} and, for $g=0$, the energy estimate
(\ref{eq:energy-est}) holds.
\end{thrm} 
One might be tempted to replace the requirement (\ref{2.1.7}) by the
weaker estimate
\begin{equation}\label{2.1.8}
\int_0^T e^{-2\eta t}\left( \|u(\cdot,t)\|^2_{\xHzero(\Omega)}+
 \|u(\cdot,t)\|^2_{\xHzero(\Gamma)}\right) dt \le \frac{K}{\eta^{\alpha}}  
\int_0^T e^{-2\eta t} \|F(\cdot,t)\|^2_{\xHzero(\Omega)}dt .
\end{equation}
However, the definition is not stable against lower order perturbations. In section \ref{sec_3-3} we
will give an example which is algebraically {\it Unstable}, i.e., with time the solution loses more
and more derivatives.

\begin{dfntn}\label{definition2.1.5}
We call the problem (\ref{2.1.1})--(\ref{2.1.3}) {\it Unstable} if the estimate (\ref{2.1.7})
does not hold.
\end{dfntn}

For first order systems the generalization to variable coefficients (and then to
quasilinear equations) uses the theory of pseudo-differential operators and
requires the construction of a symmetrizer, as described in \cite{Kreiss1970},
which is smooth in all variables. If the problem is {\it Strongly Boundary
Stable}, then, as we have mentioned above, the same construction can be used for
second order systems. If the problem is only {\it Boundary Stable}, then we have
to modify the construction.  This can be done but is technically somewhat
complicated and the details are beyond the scope of this paper. However, we will
make the result plausible.

Since the stability against lower order perturbations is crucial for the 
generalization of our results to systems with variable coefficients in general domains,
we shall give a proof in section \ref{sec_2-2}.

It is also well known that stability against lower order perturbations allows us to use
``localization'' to decompose an initial boundary value problem on a general compact
domain into a finite number of problems which are either initial values problems in the
whole space, or initial boundary value problems in the half space. We illustrate the
technique with a simple example in one dimension.

Consider the initial boundary value problem for the wave equation on the strip
\begin{equation}
u_{tt} = u_{xx}, \quad x\in[0,1], \quad t\in[0,\infty) \label{2.1.9}
\end{equation}
with initial and boundary conditions
\begin{equation}
u(x,0)=f_1(x), \quad u_t(x,0)=f_2(x), \quad B_0 u(0,t)=g_0(t), \quad B_1
u(1,t)=g_1(t), \label{2.1.10}
\end{equation}
where $B_0$ and $B_1$ are linear first order differential operators.

A partition of unity of $[0,1]$ can be chosen as a set of three $\xCinfty$
functions $\varphi_1(x),$ $\varphi_2(x),$ $\varphi_3(x)$ where $\varphi_1$ is a cutoff function
\[
\varphi_1(x)=1 \quad {\rm if }~x\le 1/4, \quad \varphi_1(x)=0 \quad {\rm if }~
x\ge 1/2.
\]
Similarily,
\[  
\varphi_3(x)=0 \quad {\rm if }~x\le 1/2, \quad \varphi_3(x)=1 \quad {\rm if }~
x\ge 3/4,
\]
and
\[
\varphi_2(x)= 1-\varphi_1(x)-\varphi_3(x).
\]
We now define the functions
\[
u_i(x,t)=\varphi_i(x) u(x,t), \quad i=1,2,3.
\]
Clearly $u(x,t)=u_1(x,t)+u_2(x,t)+u_3(x,t)$ for all $x\in[0,1]$ and
\[
u_{itt}=\varphi_i u_{tt}= \varphi_i u_{xx} = u_{ixx} + L_i,
\]
where
\[
L_i = -2\varphi_{ix} (u_{1x} + u_{2x} + u_{3x}) - \varphi_{ixx} (u_1 + u_2 + u_3),
\]
consist only of lower order terms and has the same support as $u_i.$ Defining also
$f_{1i}(x)=\varphi_i(x) f_1(x)$ and $f_{2i} = \varphi_i(x) f_2(x)$ we obtain that $u_1$
solves the half line problem
\[\begin{split}
u_{1tt} &= u_{1xx} + L_1, \quad x\in[0,\infty), \quad t\ge 0,\\ 
u_1(x,0) &= f_{11}(x), \quad u_{1t}(x,0)=f_{21}(x), \quad B_0 u_1(0,t)=g_0(t).
\end{split}
\]
$u_2$ solves the initial value problem on the whole line
\[
u_{2tt}=u_{2xx} + L_2, \quad u_2(x,0)=f_{12}(x), \quad u_{2t}(x,0)=f_{22}(x),
\quad x\in(-\infty,\infty), \quad t\ge 0,
\]
and $u_3$ solves the half line problem
\[\begin{split}
u_{3tt} &= u_{3xx} + L_3, \quad x\in(-\infty,1], \quad t\ge 0, \cr
\quad u_3(x,0) &= f_{13}(x), \quad u_{3t}(x,0)=f_{23}(x),
\quad B_1 u_3(1,t)=g_1(t).
\end{split}
\]

If the three problems for $u_1,$ $u_2$ and $u_3$ are well posed, then the
original problem (\ref{2.1.9})--(\ref{2.1.10}) is well posed.

To treat variable coefficient problems  one can
invoke what is known as the ``principle of frozen coefficients'' to replace the
problem by one with constant coefficients.  Heuristically, one can think that if
one localizes the problem to very small regions then the coefficients of the
equation in each region are nearly constant and the behavior of the solution is
near to that of an equation with constant (frozen) coefficients. The proof of
the validity of this requires the use of pseudo-differential theory. Here we
claim the validity of this principle for our problem but do not go into the
details.

%%%%%%%%%%%%%%%%%%%%%%%%%%%%%%%%%%%%%%%%%%%%%%%%%%%%%%%%%%%%%%%%%%%%%%%%
% old sec_2-2.tex
%%%%%%%%%%%%%%%%%%%%%%%%%%%%%%%%%%%%%%%%%%%%%%%%%%%%%%%%%%%%%%%%%%%%%%%%
%\setcounter{equation}{0}
\subsection{Stability against lower order perturbations}\label{sec_2-2}

If the problem is {\it Strongly Boundary Stable}, then it is {\it Strongly Stable} and therefore
stable against lower order perturbations.  We shall now prove that the corresponding results hold
for well posed problems.
\begin{thrm}\label{theorem2.2.1}
Consider the problem (\ref{2.1.1})--(\ref{2.1.3}) for $F=0,$ $g\ne 0$ and change the boundary
conditions (\ref{2.1.3}) to
\begin{equation}
Lu=g+lu,\quad |l| \hbox{ bounded.} \label{2.2.1}
\end{equation}
If the problem is {\it Boundary Stable}, then the same is true for the
perturbed problem.
\end{thrm}
\begin{proof} We consider $lu$ as part of the data. Then (\ref{2.1.5}) becomes
\begin{equation}\label{2.2.2}
\int_0^T e^{-2\eta t} \| u(\cdot,t)\|^2_{\xHzero(\Gamma)} dt \le \frac{2}{
\eta^\alpha}\Bigl(|l|^2
 \int_0^T e^{-2\eta t} \| u(\cdot,t)\|^2_{\xHzero(\Gamma)} dt+
 \int_0^T e^{-2\eta t} \| g(\cdot,t)\|^2_{\xHzero(\Gamma)} dt\Bigr).
\end{equation}
By choosing $\eta_0$ such that $\frac{K}{\eta_0^\alpha}|l|^2\le \frac{1}{4},$
we obtain an estimate also for the perturbed problem.
\end{proof}
\begin{thrm}\label{theorem2.2.2}
Consider the problem (\ref{2.1.1})--(\ref{2.1.3}) with $g=0$ and change the differential equations
to
\begin{equation}
u_{tt}=P_0(D) u+P_1(D) u+F. \label{2.2.3}
\end{equation}
Here $P_1(D)$ is a first order differential operator with bounded coefficients, i.e.,
\[
\|P_1(D)u\|^2_{\xHzero(\Omega)}\le 
K_1\|u\|^2_{\xHone(\Omega)}.
\] 
Assume that our problem is {\it Stable}. Then the perturbed problem has the same property.
\end{thrm}
\begin{proof} In the same way as in Theorem \Rref{theorem2.2.1} we consider $P_1(D)u$ as
part of the forcing and choose $\eta_0$ and $\alpha$ sufficiently large. Then the desired
estimate follows.
\end{proof}
\begin{rmrk}\label{remark:1.11}
One could argue that there are too many definitions. The reason for the concept
of boundary stabillity is that it gives the easiest test for wellposedness. If
the problem is {\it Strongly Boundary Stable} then this test provides sufficient
conditions for wellposedness. If the problem is only {\it Boundary Stable} we can use
the test to find out what kind of boundary phenomena the problem has. This
knowledge is crucial in constructing numerical methods. In the {\it Boundary Stable} case the test
provides only necessary conditions; to obtain wellposedness one needs to show
that the estimate (\ref{2.1.7}) also holds. This is always true if there is an
energy estimate for homogeneous boundary data.
\end{rmrk}
\begin{rmrk}\label{remark:1.12}
If we do not want to distinguish between {\it Strongly Stable} and {\it Stable}, we call the problem
well posed.
\end{rmrk}
%%%%%%%%%%%%%%%%%%%%%%%%%%%%%%%%%%%%%%%%%%%%%%%%%%%%%%%%%%%%%%%%%%%%%%%%
% old sec_3-1.tex
%%%%%%%%%%%%%%%%%%%%%%%%%%%%%%%%%%%%%%%%%%%%%%%%%%%%%%%%%%%%%%%%%%%%%%%%
\section{A single wave equation}\label{sec_3}
\subsection{A necessary condition for well posedness}\label{sec_3-1}
To derive necessary conditions we consider in this section the halfplane problem for the wave equation
\begin{equation}
u_{tt}=u_{xx}+u_{yy} + F(x,y,t), \quad x\ge 0,~-\infty<y<\infty,~t\ge 0
\label{3.1.1}
\end{equation}
with initial conditions, at $t=0,$
\[ u(x,y,0)=f_1(x,y),\quad u_t(x,y,0)=f_2(x,y) \]
and one of four types of boundary conditions at $x=0,~-\infty<y<\infty:$
\begin{equation}\label{3.1.2}
\begin{split}
1)& \quad u_t = au_x+bu_y+g,\quad a,b~{\rm real},~|b|<1,~a>0. \\
2)& \quad u_x=ibu_y + g, \quad b~{\rm real,}~b\neq 0,~|b|<1. \\
3)& \quad u_x = g. \\
4)& \quad u_x = bu_y+g,\quad b~{\rm real}, ~b\neq 0.\\
\end{split}
\end{equation}
The source $F$ and the data $f_j,g$ are compatible smooth functions with compact
support. We are only interested in solutions with bounded $\xLtwo$-norm and
therefore we assume
\begin{equation}
\int_{-\infty}^{\infty}\int_0^\infty |u(x,y,t)|^2~dx~dy = \|u\|^2 <\infty\quad
{\rm for~every~fixed}~t. \label{3.1.3}
\end{equation}

We start with a test to find a necessary condition such that the problem is well posed. In this
chapter and throughout the rest of the paper, $s=\eta+i\xi$ denotes a complex number where $\eta$, $\xi\in {\mathbb R}$.
\begin{lmm}\label{lemma3.1.1}
Let $F=g=0.$ The problem (\ref{3.1.1})--(\ref{3.1.2}) is not well posed if
we can find a nontrivial simple wave solution of type
\begin{equation}
u=e^{st+i\omega y}\varphi(x),\quad \|\varphi(x)\|<\infty,\quad {\rm Re\,}
s>0.\label{3.1.4}
\end{equation}
\end{lmm}
\begin{proof} If we have found such a solution, then
\[
u_\alpha= e^{s\alpha t+i\omega\alpha y} \varphi(\alpha x),\quad \alpha >0.
\]
is also a solution for any $\alpha>0.$ Since ${\rm Re\,} s >0,$ we can find 
solutions which grow arbitrarily fast exponentially.
\end{proof}

We shall now discuss whether there are such solutions. Introducing (\ref{3.1.4})
into the homogeneous differential equation (\ref{3.1.1}) and homogeneous boundary
conditions (\ref{3.1.2}) gives us
\begin{equation}
\varphi_{xx}-(s^2+\omega^2)\varphi=0,\quad \|\varphi\|<\infty. \label{3.1.5}
\end{equation}
(\ref{3.1.5}) is an ordinary differential equation with constant coefficients and
boundary conditions
\begin{equation}\label{3.1.6}
\begin{split}
1)& \quad s\varphi(0) = a\varphi_x(0) + bi\omega\varphi(0),\quad a,b~{\rm real},~|b|<1,~a>0.\cr
2)& \quad \varphi_x(0) = -b\omega\varphi(0),\quad b~{\rm real,}~b\neq 0,~|b|<1. \cr
3)& \quad \varphi_x(0) = 0.\cr
4)& \quad \varphi_x(0) = bi\omega\varphi(0),\quad b~{\rm real}, b\neq 0.\cr
\end{split}
\end{equation}

The general solution of (\ref{3.1.5}) is of the form
\begin{equation}
\varphi(x)=\sigma_1 e^{\kappa x}+\sigma_2 e^{-\kappa x},\label{3.1.7}
\end{equation}
where $\pm \kappa$ are the solutions of the characteristic equation
\[
\kappa^2-(s^2+\omega^2)=0,\quad {\rm i.e.,}\quad \kappa=\sqrt{s^2+\omega^2}.
%\eqno{()}
\]
We fix the argument of $\sqrt{}$ by
\[
-\pi<~{\rm arg}(s^2+\omega^2)\le \pi,\quad {\rm
arg}\sqrt{s^2+\omega^2}=\frac{1}{2}~{\rm arg}(s^2+\omega^2).
\]

From the general theory (also proved in Lemma \Rref{lemmaA.5} in the appendix) we know that,
there is a constant $\delta >0$ such that
\[
\re \kappa \ge \delta~\re s.
\] 
Therefore $\varphi\in \xLtwo$ if and only if $\sigma_1=0.$ Introducing (\ref{3.1.7}) into the
boundary conditions gives us
\begin{equation}\label{3.1.8}
\begin{split}
1)& \quad s=- a\kappa +i\omega b,\quad a,b~{\rm real},~|b|<1,~a>0.\cr
2)& \quad \kappa = \omega b,\quad b~{\rm real,}~b\neq 0,~|b|<1.\cr
3)& \quad \kappa =0. \cr
4)& \quad \kappa =-i\omega b,\quad b~{\rm real}, b\neq 0. \cr
\end{split}
\end{equation}
Since, by assumption, $a>0$ and $\re\kappa >0,$ there are no solutions of type
(\ref{3.1.4}) for the first kind of boundary condition. It is important to stress here
that chosing the the wrong sign for $a$ in the first type of boundary condition
results into an ill posed problem and no solution can be computed. 
The second case in (\ref{3.1.8}) implies
\[
\omega^2+s^2=\omega^2b^2,\quad {\rm i.e.,}\quad s^2=\omega^2(b^2-1).
\]
Thus there are no solutions of type (\ref{3.1.4}) for $|b|<1,~b$ real. As $\re\kappa>0$
there are no solution of type (\ref{3.1.4}) for the third boundary condition. Finally,
since $\re\kappa >0$ and $b$ is real, there are no simple wave solutions of type
(\ref{3.1.4}) for the fourth kind of boundary condition.

We have proved
\begin{thrm}\label{theorem3.1.1}
For the boundary conditions (\ref{3.1.2}) there are no simple wave
solutions of type (\ref{3.1.4}) other than the trivial solution $u \equiv 0.$
\end{thrm}

Since (\ref{3.1.5}) and (\ref{3.1.8}) define eigenvalue problems, we can phrase the theorem also as

\begin{thrm}\label{theorem3.1.1'}
The eigenvalue problems (\ref{3.1.5}) and (\ref{3.1.8}) have no eigenvalues
with ${\rm Re\,}s >0.$
\end{thrm}

We shall now introduce the concept of generalized eigenvalues. For that purpose
we write (\ref{3.1.8}) in terms of normalized variables.
\[
s'=s/\sqrt{|s|^2+\omega^2},\quad
\omega'=\omega/\sqrt{|s|^2+\omega^2},\quad
\kappa'=\kappa/\sqrt{|s|^2+\omega^2}.
\]
\begin{dfntn}\label{definition3.1.1}
Let $s'=i\xi'_0,~\omega'=\omega'_0$ be a fixed point and
consider (\ref{3.1.5}),(\ref{3.1.8}) for $s'=i\xi'_0+\eta',~\omega'=\omega'_0,~\eta'>0.$
$(i\xi'_0,\omega'_0)$ is a generalized eigenvalue for a boundary condition if in
the limit $\eta'\to 0$ the boundary condition is satisfied.
\end{dfntn}

We now calculate the generalized eigenvalues. By Lemma \ref{lemmaA.7}, there are no
generalized eigenvalues for boundary conditions of type 1).

For boundary conditions of type 2), we need to consider
\[
\lim_{\eta'\to 0} \left(\sqrt{(i\xi'_0+\eta')^2+{\omega_0'}^2} -b\omega'_0\right).
\]
As $\re \kappa'\ge 0,$ there will be a generalized eigenvalue for boundary
condition 2) if and only if $b\omega'_0 > 0$ and
\[
-{\xi'}_0^2+{\omega'}_0^2=b^2{\omega'}_0^2,\quad {\rm i.e.,} \quad \xi'_0=\pm\sqrt{1-b^2}\,\omega'_0.
\]
Since
\[
\kappa'_0=\sqrt{-{\xi'}_0^2+{\omega'}_0^2}= b\omega'_0>0,
\]
the corresponding eigenfunctions are
\[
u= e^{-|b\omega_0| x}\,e^{i\omega_0(y \pm\sqrt{1-b^2}\,t)}.
\]
They represent {\it surface waves} which decay exponentially in $x$, i.e.~in the normal
direction away from the boundary. They are important phenomena in many applications
(e.g. elastic wave equations).

For boundary conditions 3) and 4) $\xi'_0,\omega'_0$ must satisfy the relation
\[
\xi'_0=\pm\sqrt{1+b^2} \omega'_0,\quad \kappa'_0 = -i\omega'_0 b,
\]
and where the sign in the first relation is chosen so that $\xi'_0 b
\omega'_0<0$ (because $\re\kappa'\ge 0.$)
The corresponding eigenfunction is 
\begin{equation}\label{eq:eigenfcn4}
u = e^{i\omega_0(bx+y)} e^{\pm i \sqrt{1+b^2} |\omega_0| t}.
\end{equation}
For boundary condition 4) ($b\ne 0,$) they are oscillatory in $x,y,t.$ For
boundary condition 3) ($b=0$), they are constant normal to the boundary and they
are called {\it glancing waves.} They are important physical phenomena
(e.g. Maxwell's equations).

We collect the results in
\begin{thrm}\label{theorem3.1.2}
There are no generalized eigenvalues for boundary conditions of
type (1). For boundary condition (2), (3) and (4) the generalized eigenvalues are given by
\[\begin{split}
2)& \quad \xi'_0=\pm\sqrt{1-b^2}\,\omega'_0,\quad \kappa'_0=b\omega'_0>0.\\
3)& \quad |\xi'_0|= |\omega'_0|, \quad \kappa'_0 = 0. \\
4)& \quad |\xi'_0|= \sqrt{1+b^2}\,|\omega'_0|,\quad \kappa'_0=-i\omega'_0 b,\quad
\xi'_0 b \omega'_0 < 0.
\end{split}
\]
\end{thrm}

%%%%%%%%%%%%%%%%%%%%%%%%%%%%%%%%%%%%%%%%%%%%%%%%%%%%%%%%%%%%%%%%%%%%%%%%
% old sec_3-2.tex
%%%%%%%%%%%%%%%%%%%%%%%%%%%%%%%%%%%%%%%%%%%%%%%%%%%%%%%%%%%%%%%%%%%%%%%%
\subsection{Reduction to a first order system of pseudo-differential equations}\label{sec_3-2} 

The estimates obtained in this and subsequent sections are expressed in Fourier-Laplace
transformed space. It is clear that all these estimates have their counterpart in physical
space (such as the estimates in the definitions of section \ref{sec_2-1}). To understand the relation
between both types of estimates we refer to chapter 7.4 of \cite{Kreiss-Lorenz1989} and
chapter 10 of~\cite{Gustafsson-Kreiss-Oliger1995}.

We consider (\ref{3.1.1})--(\ref{3.1.3}) with homogeneous initial data. We Laplace
transform the problem with respect to $t,$ Fourier transform it with respect to $y,$ and
denote the dual variables by $s,\omega$, respectively. For Re$\,s>0$ we obtain
\begin{equation}
\hat u_{xx}= (s^2 + \omega^2)\hat u + \hat F, \quad \hat u=\hat
u(x,\omega,s),~\hat F=\hat F(x,\omega,s),~ 0\le x< \infty, \label{3.2.1}
\end{equation}
with one of the boundary conditions
\begin{equation}\label{3.2.2}
\begin{split}
1)&\quad s\hat u=a\hat u_x+ib\omega\hat u+\hat g,\cr
2)&\quad \hat u_x=-b\omega \hat u+\hat g,\cr
3)&\quad \hat u_x=\hat g,\cr
4)&\quad \hat u_x=ib\omega\hat u+\hat g,\cr
\end{split}
\end{equation}
and $\|\hat u(\cdot,\omega,s)\|^2<\infty.$

Introducing a new variable by
\begin{equation}
\hat u_x=\sqrt{|s|^2+\omega^2}\,\hat v, \label{3.2.3}
\end{equation}
we write the Fourier and Laplace transformed system as a first
order system
\begin{equation}
\hat{\bf u}_x= \sqrt{|s|^2+\omega^2}\,M\hat{\bf u} + \hat{\bf F}_0.\label{3.2.4}
\end{equation}
Here
\[
\hat{\bf u}=\begin{pmatrix}\hat u\cr \hat v\end{pmatrix},\quad M=\begin{pmatrix}0 & 1\cr {\kappa'}^2 &
0\end{pmatrix},\quad \hat{\bf F}_0=
\frac{1}{\sqrt{|s|^2+\omega^2}}\begin{pmatrix}0\cr \hat F\end{pmatrix},
\]
with
\[
\kappa'=\sqrt{(s')^2+(\omega')^2}, \quad
s'=\frac{s}{\sqrt{|s|^2+\omega^2}},\quad \omega'= \frac{\omega}{\sqrt{|s|^2+\omega^2}}.
\]
The eigenvalues $\mu$ of $M$ are $\mu_1=-\kappa'$ and $\mu_2=\kappa'.$ The
boundary conditions at $x=0$ become
\begin{equation}\label{3.2.5}
\begin{split}
1)& \quad s'\hat u=a\hat v+ib\omega'\hat u+g', \quad a>0,~|b|<1,\cr
2)&\quad \hat v=-b\omega'\hat u + g', \quad |b|<1,~b\ne 0, \cr
3)& \quad \hat v=g', \cr
4)& \quad \hat v=ib\omega'\hat u+g', \quad b~\mbox{real,}~b\neq 0,\cr
\end{split}
\end{equation}
with $g'= \hat g/\sqrt{|s|^2+\omega^2}.$

\begin{rmrk}\label{symmetrizer_remark}
We present in this and the following section the easiest way to obtain the
estimates at the boundary and in the interior of the domain. To generalize these
results to variable coefficients pseudo-differential theory is needed. The
transformations $S$ and $T$ introduced below (see eqns.
(\ref{3.2.7},\ref{3.3.1})) need to be smooth in the dual variables. The
smoothness condition may fail only at the double root of $M.$ In this case the
Kreiss' symmetrizer is used to get the estimates as explained in
\cite{Kreiss1970}.
\end{rmrk}

We shall now calculate the solution for the case when $F=0,$ and
estimate it on the boundary.

The eigenvector of $M$ connected with $-\kappa'$ is given by
\begin{equation}
\bf x=\begin{pmatrix}1 \cr -\kappa'\end{pmatrix}. \label{3.2.6}
\end{equation}
The transformation
\begin{equation}
S=\begin{pmatrix}1& \bar\kappa'\cr -\kappa' & 1\end{pmatrix} \label{3.2.7}
\end{equation}
is, except for a trivial normalization, unitary and transforms $M$ 
into upper triangular form, i.e.,
\begin{equation}
S^{-1}MS=\begin{pmatrix}-\kappa' &d\cr 0 & \kappa'\end{pmatrix}, \quad
d=\frac{1+|\kappa'|^4}{1+|\kappa'|^2}. \label{3.2.8}
\end{equation}
Here $S,S^{-1},d$ are uniformly bounded and depend smoothly on $\kappa'.$

Introducing a new variable by
\begin{equation}
\begin{pmatrix}\hat u\cr \hat v\end{pmatrix}= S \begin{pmatrix}\tilde u\cr
\tilde v\end{pmatrix} \label{3.2.9}
\end{equation}
transforms (\ref{3.2.4}) into 
\begin{equation}
\begin{pmatrix}\tilde u\cr \tilde v\end{pmatrix}_x= \begin{pmatrix}-\kappa
&d\sqrt{|s|^2+\omega^2}\cr 0 & \kappa\end{pmatrix}
 \begin{pmatrix}\tilde u\cr \tilde v\end{pmatrix}. \label{3.2.10}
\end{equation}
As the solution is in $\xLtwo,$ we have $\tilde v=0$ and also $\hat u=\tilde u.$
Thus the boundary conditions become
\begin{equation}\label{3.2.11}
\begin{split}
1)&\quad (s'+a\kappa'-ib\omega')\tilde u=g',\\
2)&\quad (\kappa'-b\omega')\tilde u=-g',\\
3)&\quad \kappa'\tilde u=-g',\\
4)&\quad (\kappa'+ib\omega')\tilde u=-g'.\\
\end{split}
\end{equation}
By (\ref{3.1.8}), these boundary conditions become singular exactly at the generalized
eigenvalues.

\begin{rmrk}\label{remark3.2.1}
From all the boundary conditons the first one is the most benign. By Lemma
\ref{lemmaA.7} we obtain the estimate on the boundary
\begin{equation}
|\tilde u(0,\omega,s)|^2 \le \frac{{\rm const.}}{|s|^2 + \omega^2} |\hat
g(\omega,s)|^2. \label{3.2.12}
\end{equation}
In this case we gain a derivative on the boundary and the problem is {\it Strongly Boundary
  Stable}. According to the classical theory \cite{Kreiss1970,Kreiss-Lorenz1989} the problem is {\it
  Strongly Stable}.  Moreover, the principle of localization holds and the problem can be
generalized to variable coefficients and then to quasilinear equations. It is worth noticing here
that away from generalized eigenvalues, i.e. when the coefficients on the left hand side of
(\ref{3.2.11}) are strictly away from zero, the estimate (\ref{3.2.12}) holds also for boundary
conditions 2), 3) and 4) and therefore the problem can be treated by the classical theory.
\end{rmrk}

Because of the previous remark, we only need to study the estimates near the
generalized eigenvalues. We have

\begin{thrm}\label{theorem3.2.1}
The problem (\ref{3.1.1})--(\ref{3.1.3}) with $F=0$ and $f_1=f_2=0~$ has a
unique solution in $\xLtwo$ which, Fourier-Laplace transformed, is given by
\begin{equation}
\hat u = e^{-\kappa x}\hat u(0,\omega,s),\quad
\kappa=\sqrt{s^2+\omega^2}.\label{3.2.13}
\end{equation}
For the different boundary conditions sharp estimates follow. For boundary
condition 1) and, ``away'' from generalized eigenvalues, for all other boundary
conditions the problem is {\it Strongly Boundary Stable} and
\begin{equation}
1)\quad |\hat u(0,\omega,s)|^2\le \frac{{\rm const.}}{|s|^2+\omega^2}|\hat g|^2.
\label{3.2.14}
\end{equation}
Near generalized eigenvalues the estimates for boundary conditions 2), 3) and 4) are
\begin{equation}\label{3.2.15}
\begin{split}
2)&\quad |\hat u(0,\omega,s)|^2\le \frac{\rm const.}{\eta^2} |\hat g|^2, \\
3)&\quad |\hat u(0,\omega,s)|^2\le {\rm const.}\,\frac{|\hat g|^2}{|\kappa|^2}\le
{\rm const.} \frac{|\hat g|^2}{\eta(|s|^2+\omega^2)^{1/2}},\\
4)&\quad |\hat u(0,\omega,s)|^2\le\frac{{\rm const.}}{\eta^2} |\hat g|^2,
\end{split}
\end{equation}
and the problem is {\it Boundary Stable}.
\end{thrm}

\begin{proof} Clearly, as $\tilde v=0$ and $\hat u=\tilde u,$
(\ref{3.2.13}) is the only solution to (\ref{3.2.10}). We need to consider only a neighbourhood of the
generalized eigenvalues $(i\xi'_0,\omega'_0).$ 
For the second boundary condition $\xi'_0=\pm\sqrt{1-b^2}\,\omega'_0,$ i.e.,
\[
s'=i(\xi'_0+\tilde{\xi}')+\eta',\quad \eta'\ge 0,\quad \omega'=\omega'_0+
{\tilde \omega}',\quad |{\tilde\xi}'|+|\tilde\omega'|+\eta' \ll 1.
\]
Since $\kappa'\ne 0$ at $(i\xi'_0,\omega'_0),$ we can use Taylor expansion.
A simple perturbation caculation shows that the worst estimate occurs
for $\tilde\xi'=\tilde\omega'=0.$ In this case we have,
\begin{equation}\label{3.2.16}
\begin{split}
|\kappa'-b\omega'|&=\Bigl|\sqrt{-\xi^{\prime 2}_0+2i\xi'_0\eta'+\eta^{\prime 2}+\omega^{\prime 2}_0}
-b\omega'_0\Bigr|\\
&=\Bigl|\sqrt{b^2\omega^{\prime 2}_0+2i\xi'_0\eta'+\eta^{\prime 2}} -b\omega'_0\Bigr|\\
&\approx \bigg||b\omega'_0|-b\omega'_0+ \frac{i\xi'_0\eta'}{|b\omega'_0|}\bigg|
\approx \frac{\sqrt{1-b^2}}{|b|}\,\eta'\quad {\rm if}\quad b\omega'_0>0.
\end{split}
\end{equation}
Thus we have
\[
|\tilde u(0,\omega,s)|^2\le{\rm const.} \frac{|g'|^2}{{\eta'}^2}
={\rm const.} \frac{|\hat g|^2}{\eta^2}.
\]
A similar perturbation calculation gives
\begin{equation}
|\kappa' + ib\omega'| \simeq \frac{\sqrt{1+b^2}}{|b|}\eta' + {\cal
O}(\eta'^2) \label{3.2.17}
\end{equation}
which gives, for the last boundary condition
\[
|\tilde u(0,\omega,s)|^2\le{\rm const.} \frac{|\hat g|^2}{\eta^2}.
\]
For the third boundary condition on can do better. Lemma A.4 with $b=0,$ gives us
\begin{equation}\label{3.2.18}
|\tilde u(0,\omega,s)|^2\le{\rm const.}\left(\frac{|g'|}{|\kappa'|}\right)^2
={\rm const.}\left(\frac{|\hat g|}{|\kappa|}\right)^2
\le \frac{\rm const.}{\eta}\,\frac{|\hat g|^2}{(|s|^2+\omega^2)^{1/2}}.
\end{equation}
This proves the theorem.
\end{proof}

\begin{thrm}\label{theorem3.2.2}
When $f_1 = f_2 = 0$ and $F=0.$ The unique
solution to our problem, described in Theorem \Rref{theorem3.2.1}, satisfies the
following interior estimates. For boundary condition 1) and, away from
generalized eigenvalues, for boundary conditions 2), 3) and 4)
\begin{equation}
\|\tilde u\|^2 \le \frac{{\rm const.}}{\eta}~\frac{|\hat
g|^2}{|s|^2+\omega^2}.\label{3.2.19}
\end{equation}
Close to generalized eigenvalues, for the corresponding boundary conditions, we
have
\begin{equation}\label{3.2.20}
\begin{split}
2)&\quad\|\tilde u\|^2 \le \frac{{\rm const.}}{\eta^2} \frac{|\hat
g|^2}{(|s|^2+\omega^2)^{1/2}}\\
3)&\quad\|\tilde u\|^2 \le \frac{{\rm const.}}{\eta^{3/2}} \frac{|\hat g|^2}{(|s|^2+\omega^2)^{3/4}}.\\
4)&\quad\|\tilde u\|^2 \le \frac{{\rm const.}}{\eta^3} |\hat g|^2
\end{split}
\end{equation}
\end{thrm}

\begin{proof} By Theorem \Rref{theorem3.2.1} and Lemma A.2, the solution satifies
\[
\|\tilde u\|^2\le \frac{1}{\re\kappa} |\tilde u(0,\omega,s)|^2.
\]
Since, by Lemma A.5, always $\re \kappa\ge \delta_4\eta,$ we obtain from
(\ref{3.2.13}) that (\ref{3.2.19}) is valid. 2) follows because, by Theorem
\Rref{theorem3.1.2},
$\re\kappa'\simeq b\omega'_0$ and then $\re\kappa\ge{\rm const.}\,\sqrt{|s|^2+\omega^2}.$ 
 3) follows from (\ref{3.2.18}), Lemma A.6 and Lemma A.4, according to
\begin{equation}\label{3.2.21}
\|\tilde u\|^2\le {\rm const.} \frac{|\hat g|^2}{\re\kappa|\kappa|^2}\le
{\rm const.}\frac{|\hat g|^2}{\eta |\kappa| \sqrt{|s|^2+\omega^2}}\le {\rm
const.} \frac{|\hat g|^2}{\eta^{3/2}(|s|^2+\omega^2)^{3/4}}.
\end{equation}
Finally 4) corresponds to 4) of (\ref{3.2.15}). This proves the theorem.
\end{proof}

%%%%%%%%%%%%%%%%%%%%%%%%%%%%%%%%%%%%%%%%%%%%%%%%%%%%%%%%%%%%%%%%%%%%%%%%
% old sec_3-3.tex
%%%%%%%%%%%%%%%%%%%%%%%%%%%%%%%%%%%%%%%%%%%%%%%%%%%%%%%%%%%%%%%%%%%%%%%%
\subsection{Estimates for homogeneous boundary data}\label{sec_3-3}
We consider now the problem (\ref{3.2.4}),(\ref{3.2.5}) with $g'=0$ and treat
only the cases 2), 3) and 4) where there are generalized eigenvalues (see Remark
\ref{remark3.2.1}).

For $\eta=\re s >0,$ the eigenvalues of $M$ are distinct and therefore we can 
transform (\ref{3.2.4}) to diagonal form by the transformation
\begin{equation}\label{3.3.1}
T=\begin{pmatrix}1 & 1 \cr -\kappa' & \kappa'\end{pmatrix}, \quad T^{-1}=\frac{1}{2}\begin{pmatrix}1
& -1/\kappa' \cr 1 & +1/\kappa'\end{pmatrix},
\end{equation}
Let
$\tilde u, \tilde v$ be defined by
\begin{equation}
\begin{pmatrix}\hat u\cr \hat v\end{pmatrix} = T\begin{pmatrix}\tilde u\cr
\tilde v\end{pmatrix}. \label{3.3.2}
\end{equation}
Then, (\ref{3.2.4}) becomes
\begin{equation}
\begin{pmatrix}\tilde u\cr \tilde v\end{pmatrix}_x=\begin{pmatrix}-\kappa & 0\cr
0 & \kappa\end{pmatrix}\begin{pmatrix}\tilde
u\cr \tilde v\end{pmatrix} + \frac{1}{
2\kappa'\sqrt{|s|^2+\omega^2}}\begin{pmatrix}-\hat F\cr 
\hat F\end{pmatrix},\label{3.3.3}
\end{equation}
with boundary conditions 
\begin{equation}\label{3.3.4}
\begin{split}
2)&\quad (\kappa'-b\omega')\tilde u(0,\omega,s)=(\kappa'+b\omega')\tilde
v(0,\omega,s).\\
3)&\quad \tilde u(0,\omega,s)=\tilde v(0,\omega,s).\\
4)&\quad (\kappa'+ib\omega')\tilde u(0,\omega,s)=(\kappa'-ib\omega')\tilde
v(0,\omega,s).
\end{split}
\end{equation}

The equations (\ref{3.3.3}) are decoupled and as $\re\kappa>0,$ Lemma A.1 gives 
for all boundary conditions
\begin{equation}\label{3.3.5}
\begin{split}
|\tilde v(0,\omega,s)|^2&\le \frac{1}{8|\kappa'|^2\re\kappa} \frac{\|\hat F\|^2}{
|s|^2+\omega^2} = \frac{1}{8|\kappa'|^2\re\kappa'} \frac{\|\hat F\|^2}{
(|s|^2+\omega^2)^{3/2}} ,\\
 \|\tilde v(\cdot,\omega,s)\|^2&\le \frac{1}{4|\kappa'|^2|\re\kappa|^2} \frac{\|\hat
F\|^2|}{|s|^2+\omega^2}=\frac{1}{4|\kappa'|^2|\re\kappa'|^2} \frac{\|\hat F\|^2}{(|s|^2+\omega^2)^2}.
\end{split}
\end{equation}

We use the boundary conditions to estimate $\tilde u(0,\omega,s).$

Theorem \Rref{theorem3.1.2} tells us that $\re\kappa'_0\approx b\omega'_0 >0$ in a neighborhood of
the generalized eigenvalue connected with boundary condition 2). Therefore the 
perturbation calculation (\ref{3.2.16}) gives us
\begin{equation}\label{3.3.6}
\begin{split}
|\tilde u(0,\omega,s)|^2&=\Bigl|\frac{\kappa'+b\omega'}{
\kappa'-b\omega'}\Bigr|^2\,|\tilde v(0,\omega,s)|^2
\le \frac{{\rm const.}}{\eta^{\prime 2}}\frac{\|\hat F\|^2}{(|s|^2+\omega^2)^{3/2}} \\
&\le \frac{{\rm const.}}{\eta^{2}}\frac{\|\hat F\|^2}{(|s|^2+\omega^2)^{1/2}}.
\end{split}
\end{equation}

The interior estimate in this case follows from Lemma A.2, $\re
\kappa=\re\kappa'\sqrt{|s|^2+\omega^2}$ and (\ref{3.3.6})
\begin{equation}\label{3.3.7}
\begin{split}
\|\tilde u(\cdot,\omega,s)\|^2&\le \frac{|\tilde u(0,\omega,s)|^2}{\re\kappa}+
\frac{\rm const.}{(\re\kappa)^2} \frac{\|\hat F\|^2}{ |s|^2+\omega^2}\\
&\le \frac{\rm const.}{ \eta^2} \frac{\|\hat F\|^2}{|s|^2+\omega^2}
\end{split}
\end{equation}
Since the transformation $T$ is bounded, (\ref{3.3.2}) tells us that the estimates
(\ref{3.3.5})--(\ref{3.3.7}) are also valid for $\hat u,\hat v.$ Thus the problem is {\it Stable}.

For boundary condition 3), the generalized eigenvalue is $\kappa'=0$ and therefore,
for $\re s'>0,$ we know only that $\re\kappa'\ge\delta\re s'.$ However, by Lemma A.6, we have a strong
estimate for $|\kappa|\re\kappa$ and the estimates (\ref{3.3.5}) become
\begin{equation}\label{3.3.8}
\begin{split}
|\tilde v(0,\omega,s)|^2&\le \frac{{\rm const.}}{ \eta|\kappa|}\,\frac{\|\hat
F\|^2}{ \sqrt{|s|^2+\omega^2}}
\le \frac{{\rm const.}}{\eta^2}
\,\frac{\|\hat F\|^2}{ \sqrt{|s|^2+\omega^2}}\\
\|\tilde v(\cdot,\omega,s)\|^2&\le 
 \frac{{\rm const.}}{\eta^2}
\,\frac{\|\hat F\|^2}{ |s|^2+\omega^2}.
\end{split}
\end{equation}
By (\ref{3.3.4}), the same estimate holds for $\tilde u(0,\omega,s).$ By Lemma
A.2, we obtain the interior estimate
\begin{equation}
\|\tilde u(\cdot,\omega,s)\|^2 \le 
 \frac{{\rm const.}}{\eta^2} \,\frac{\|\hat F\|^2}{ |s|^2+\omega^2}.\label{3.3.9}
\end{equation}
Since $\hat u,\hat v$ satisfy the same estimates, the problem is {\it Stable}.

For boundary condition 4), the generalized eigenvalue is $\kappa'_0=-i\omega'_0 b,~b\ne 0,$
is purely imaginary. Therefore, for $\re s'>0,$ we can only use the estimate 
$\re\kappa'\ge\delta \eta'.$ Instead of (\ref{3.3.6}) we obtain now
\[
|\tilde u(0,\omega,s)|^2=\Bigl|\frac{\kappa'-i\omega' b}{ \kappa'+i\omega' b}\Bigr|^2\,
|\tilde v(0,\omega,s)|^2\le 
 \frac{{\rm const.}}{\eta^3} \|\hat F\|^2.
\]
Again, the same estimates hold for $\hat u,\hat v.$ The estimates are sharp. Therefore
we do not obtain the desired interior estimate and the problem is {\it Unstable}. We have proved
\begin{thrm}\label{theorem3.3.1}
For boundary conditions 2) and 3) our problem is {\it Stable}, but not for boundary condition 4).
\end{thrm}

\begin{rmrk}\label{remark3.3.1}
The estimates obtained for boundary conditions 2) and 3) tell us that
$\hat u$ gains one derivative with respect to the forcing in the interior of the
domain and half a derivative on the boundary. The problem can be localized and
generalized to variable coefficients and then to quasilinear equations. On the
other hand, the estimates for the problem with boundary condition 4) show that
not even a fractional derivative is gained with respect to the forcing. This,
for a second order equation, means that one derivative of the solution is lost at
every reflection on the boundary. The problem can not be localized. We do not
pursue this problem but illustrate below this bad type of behavior with a
simple example: a first order system with a boundary condition equivalent to 4).
\end{rmrk}

\noindent{\bf An example:} Boundary reflection with loss or gain of differentiability. Consider a
system of differential equations
\begin{equation} 
u_t=-u_x,\quad v_t=v_x,\quad 0\le x\le 1,\quad t>0, \label{3.3.17}
\end{equation}
with boundary conditions
\begin{equation} 
u(0,t)=v_x(0,t),\quad v(1,t)=u_x(1,t). \label{3.3.18}
\end{equation}
Then
\begin{equation} 
u=e^{\lambda(t-x)}u_0,\quad v=e^{\lambda(t+x)} v_0, \label{3.3.19} 
\end{equation}
is a solution of (\ref{3.3.17}).

Introducing (\ref{3.3.19}) into (\ref{3.3.18}) gives us
\begin{equation}
u_0=\lambda v_0,\quad e^\lambda v_0=-\lambda e^{-\lambda} u_0.
\label{3.3.20} 
\end{equation}
Thus we obtain a solution of (\ref{3.3.17}),(\ref{3.3.18}) if
\begin{equation}
e^{2\lambda}=-\lambda^2. \label{3.3.21} 
\end{equation}
Let
\[ \lambda=\lambda_n=\pi in+\tilde\lambda_n,\quad n=1,2,\ldots,\]
then (\ref{3.3.21}) becomes
\begin{equation}
e^{2\tilde \lambda_n}=\pi^2 n^2-2\pi in\tilde \lambda_n-\tilde\lambda_n^2.
\label{3.3.22} 
\end{equation}
(\ref{3.3.22}) has a solution 
\[ \tilde\lambda_n\approx {\rm log}\,\pi n.\]
Therefore the solution (\ref{3.3.19}) grows like
\begin{equation}
e^{\lambda t}\simeq e^{\pi i nt}\cdot e^{t {\rm log}\,\pi n}
= e^{\pi int} (\pi n)^t. \label{3.3.23}
\end{equation}
If the initial data can be expanded into a Fourier series
\[ u(x,0)=\sum_n e^{\lambda_nx}\hat u(\lambda_n), \]
then (\ref{3.3.23}) tells us that the solution loses more and more derivatives with time.

Now change the boundary conditions (\ref{3.3.18}) to
\[ u_x(0,t)=v(0,t),\quad  v_x(1,t)=u(1,t). \]
Then we obtain, instead of (\ref{3.3.21}),
\[
-\lambda u_0=v_0,\quad \lambda e^{\lambda} v_0=e^{-\lambda}u_0,
\]
i.e.,
\begin{equation}
e^{2\lambda}=-\frac{1}{\lambda^2}.  \label{3.3.24}
\end{equation}
Therefore there is no loss of derivatives.

Geometrically, the two sets of boundary conditions represents two
different situations. In the first case any wave loses a derivative
when reflected at the boundary. In the second case, it gains a derivative.

%%%%%%%%%%%%%%%%%%%%%%%%%%%%%%%%%%%%%%%%%%%%%%%%%%%%%%%%%%%%%%%%%%%%%%%%
% old sec_4-1.tex
%%%%%%%%%%%%%%%%%%%%%%%%%%%%%%%%%%%%%%%%%%%%%%%%%%%%%%%%%%%%%%%%%%%%%%%%
%\newpage
%\setcounter{equation}{0}
\section{Second order systems of hyperbolic equations}\label{sec_4}
\subsection{Two wave equations}\label{sec_4-1}
In this section we consider two wave equations coupled through the boundary
conditions.
\begin{equation}
u_{1tt} = u_{1xx} + u_{1yy}, \qquad u_{2tt} = u_{2xx} + u_{2yy}, \label{4.1.1}
\end{equation}
on the halfplane $x\ge 0,$ $-\infty < y < \infty,$ for $t\ge 0,$ with
homogeneous initial conditions
\begin{equation}
u_i(x,y,0)=u_{it}(x,y,0)=0, \quad i=1,2, \quad t=0, \label{4.1.2}
\end{equation}
and boundary conditions at $x=0,$
\begin{equation}
u_{1x} + b_1 u_{2y} = g_1, \quad u_{2x} + b_2 u_{1y}=g_2, \quad x=0.
\label{4.1.3}
\end{equation}
Here $b_1,$ $b_2,$ are real and $g_j=g_j(y,t),$~$j=1,2,$ are smooth
functions which are compatible with the initial data (for example, functions that
vanish near $t=0$). We want to show that our techniques of section \ref{sec_3}
can still be used to describe the behavior of the solution.

Fourier and Laplace transform lead to
\begin{equation}
\hat u_{1xx} - (\omega^2 + s^2)\hat u_1=0,\quad \hat u_{2xx} - (\omega^2 +
s^2)\hat u_2=0, \quad \re s>0. \label{4.1.4}
\end{equation}
Thus we obtain solutions that belong to $\xLtwo$
\begin{equation}
\hat u_1=e^{st+i\omega y-\kappa x} u_{10},\quad \hat u_2=e^{st+i\omega y-\kappa
x} u_{20}, \label{4.1.5}
\end{equation}
where
\[
\kappa = \sqrt{\omega^2 + s^2}, \quad {\rm for }~\re s>0.
\]
We recall here that $\re \kappa > 0$ when $\re s>0.$ The transformed boundary
conditions become
\begin{equation}\label{4.1.7}
\begin{split}
-\kappa u_{10} + i\omega b_1 u_{20} &= \hat g_1,\\
i\omega b_2 u_{10} -\kappa u_{20} &= \hat g_2.
\end{split}
\end{equation}
A simple calculation shows that (\ref{4.1.7}) has a unique solution
\begin{equation}
u_{10} = -\frac{\kappa\hat g_1 + i\omega b_1\hat g_2}{ \kappa^2 + \omega^2
b_1b_2}, \quad u_{20} = -\frac{\kappa\hat g_2 + i\omega b_2\hat g_1}{ \kappa^2 + \omega^2
b_1b_2},\label{4.1.8}
\end{equation}
if and only if
\begin{equation}
{\rm Det}\begin{pmatrix}-\kappa & i\omega b_1\cr i\omega b_2 &
-\kappa\end{pmatrix} = \kappa^2 +
\omega^2 b_1b_2 \neq 0.\label{4.1.9}
\end{equation}
There is an eigenvalue to the homogeneous problem (\ref{4.1.4}), (\ref{4.1.7}) with $\hat
g_1=\hat g_2=0,$ if the homogeneous system (\ref{4.1.7}) has a nontrivial solution. By
(\ref{4.1.7}), this is the case if
\begin{equation}
\kappa^2 + \omega^2 b_1b_2 = s^2 + \omega^2 (b_1b_2 + 1)=0.\label{4.1.10}
\end{equation}
There are five different situations:
\begin{enumerate}
\item $b_1b_2 < -1.$ By (\ref{4.1.10}), there are eigenvalues $s$ with $\re s>0.$
Therefore our problem is not well posed.

\item $b_1b_2 = -1.$ Now $s=0$ is a generalized eigenvalue. For $s\neq 0$ the
solution (\ref{4.1.8}) becomes
\begin{equation}
u_{10} = -\frac{\kappa\hat g_1 + i\omega b_1\hat g_2}{ s^2}, \quad u_{20} =
-\frac{\kappa\hat g_2 + i\omega b_2\hat g_1}{ s^2}.\label{4.1.11}
\end{equation}
In general, we can not expect that $\hat u_{10},~\hat u_{20}$ stay bounded for
$s\to 0.$ We need to assume that $g_j=\tilde g_{jtt}$ are the second time
derivatives of smooth functions.

\item $-1<b_1b_2<0.$ By (\ref{4.1.5}),(\ref{4.1.10}), we obtain generalized
eigenvalues and eigenfunctions if $s=\pm i \sqrt{1+b_1b_2}\,|\omega|,$
$\kappa=\sqrt{|b_1b_2|}\,\omega.$ The solution (\ref{4.1.8}) is singular at the
eigenvalues and $\hat u_{10},~\hat u_{20}$ have a first order pole. In physical
space we obtain surface waves whose amplitudes become large for $b_1b_2 \to
-1.$

\item $b_1b_2=0.$ Now $s=\pm \omega,$ $\kappa=0$ determine the generalized
eigenvalues and eigenfunctions. The behavior is the same as for the Neuman
problem. We obtain glancing waves.

\item $b_1b_2>0.$ In this case the generalized eigenvalues and eigenfunctions
are determined by $s=\pm i \sqrt{1+b_1b_2}\,\omega,$ $\kappa=\pm i
\sqrt{b_1b_2}\,\omega.$ The eigenfunctions are oscillatory and the solution
behaves like the solution of the wave equation in section \ref{sec_3-3} with
boundary condition 4). Thus it is {\it Unstable}.
\end{enumerate}
In summary we can state that the problem (\ref{4.1.1})--(\ref{4.1.3}) is well
posed for boundary conditions 2)--4) and {\it Unstable} for boundary condition
5). Also, there can be numerical diffculties if $b_1 b_2 +1$ is zero or close to
zero.

When we study the estimates for the problem in the cases 3) and 4) we get
completely analogous estimates as the ones found for a single wave equation. 

Consider the system
\begin{equation}
u_{1tt}=u_{1xx}+u_{1yy}-F_1,\quad u_{2tt}=u_{2xx}+u_{2yy}-F_2, \label{4.1.12}
\end{equation}
with boundary conditions
\begin{equation}
u_{1x}+b_1 u_{2y}=0,\quad u_{2x}+b_2 u_{1y}=0, \quad x=0. \label{4.1.13}
\end{equation}
Fourier and Laplace transform gives us
\[
\hat u_{1xx}=(s^2+\omega^2)\hat u_1 +\hat F_1,\quad 
\hat u_{2xx}=(s^2+\omega^2)\hat u_2 +\hat F_2,
\]
and define 
\[
\hat u_{1x}=\sqrt{|s|^2+\omega^2}\, \hat v_1,\quad
\hat u_{2x}=\sqrt{|s|^2+\omega^2}\, \hat v_2.
\] 
Thus, the first order form of the problem becomes
\begin{equation}\label{4.1.14}
\begin{split}
\begin{pmatrix}\hat u_1 \\ \hat v_1\end{pmatrix}_x&=\sqrt{|s|^2+\omega^2}\, 
\begin{pmatrix} 0 & 1\cr \kappa'^2 & 0\end{pmatrix}
\begin{pmatrix}\hat u_1\cr \hat v_1\end{pmatrix}+\frac{1}{\sqrt{|s|^2+\omega^2}}
\,\begin{pmatrix}0\cr \hat F_1\end{pmatrix}\cr 
\begin{pmatrix}\hat u_2\cr \hat v_2\end{pmatrix}_x&=\sqrt{|s|^2+\omega^2}\, 
\begin{pmatrix} 0 & 1\cr \kappa'^2 & 0\end{pmatrix}
\begin{pmatrix}\hat u_2\cr \hat v_2\end{pmatrix}+\frac{1}{\sqrt{|s|^2+\omega^2}}
\,\begin{pmatrix}0\cr \hat F_2\end{pmatrix}
\end{split} 
\end{equation}
with boundary conditions
\[\begin{split}
\sqrt{|s|^2+\omega^2}\,\hat v_{10} + i\omega b_1\hat u_{20} &=0, \\
\sqrt{|s|^2+\omega^2}\,\hat v_{20} + i\omega b_2\hat u_{10} &=0.
\end{split}
\]
We transform (\ref{4.1.14}) to diagonal form. Let
\[
\begin{pmatrix}\hat u_j\cr \hat v_j\end{pmatrix}=\begin{pmatrix}1 & 1\cr
-\kappa'& \kappa'\end{pmatrix} \begin{pmatrix}\tilde u_j\cr \tilde
v_j\end{pmatrix},\quad j=1,2.
\]
Then
\begin{equation}\label{4.1.15}
\begin{split}
 \begin{pmatrix}\tilde u_1\cr \tilde v_1\end{pmatrix}_x
 &=\begin{pmatrix}-\kappa &0\cr 0 & \kappa\end{pmatrix}
 \begin{pmatrix}\tilde u_1\cr \tilde v_1\end{pmatrix}_x
+\frac{1}{\sqrt{|s|^2+\omega^2}} \tilde F_1\cr
 \begin{pmatrix}\tilde u_2\cr \tilde v_2\end{pmatrix}_x
 &=\begin{pmatrix}-\kappa &0\cr 0 & \kappa\end{pmatrix}
 \begin{pmatrix}\tilde u_2\cr \tilde v_2\end{pmatrix}_x
+\frac{1}{\sqrt{|s|^2+\omega^2}} \tilde F_2\cr
\end{split}
\end{equation}
with boundary conditions
\begin{equation}\label{4.1.16}
\begin{split}
-\kappa\tilde u_{10}+ i \omega b_1\tilde u_{20}&=-\kappa \tilde v_{10},\\
i \omega b_2\tilde u_{10}-\kappa \tilde u_{20} &=-\kappa \tilde v_{20}.
\end{split}
\end{equation}

In the same way as in section \ref{sec_3-3} we can determine $-\kappa\tilde v_{10},$
$-\kappa \tilde v_{20}$ by solving the equations (\ref{4.1.15}) for $\tilde v_1,
\tilde v_2$ and reduce the problem (\ref{4.1.12}),(\ref{4.1.13}) to the previous problem
(\ref{4.1.1})--(\ref{4.1.3}). Thus we obtain the same estimates.

%%%%%%%%%%%%%%%%%%%%%%%%%%%%%%%%%%%%%%%%%%%%%%%%%%%%%%%%%%%%%%%%%%%%%%%%
% old sec_5-1.tex
%%%%%%%%%%%%%%%%%%%%%%%%%%%%%%%%%%%%%%%%%%%%%%%%%%%%%%%%%%%%%%%%%%%%%%%%
\subsection{General systems. The Cauchy problem}\label{sec_5-1}
We consider the Cauchy problem for the homogeneous system of the introduction
\begin{equation}
u_{tt}=P_0(D) u,\quad t\ge 0,~x\in \xR^r:~-\infty<x_j<\infty,~j=1,2,\dots r
\label{5.1.1}
\end{equation}
were
\begin{equation}
P_0(D)= A_1D_1^2+\sum_{j=2}^r B_jD_j^2. \label{5.1.2}
\end{equation}
$A_1=A_1^* >0,$~$B_j=B_j^*>0,$ are positive definite symmetric $n\times n$
matrices and $u$ is a vector valued function with $n$ components.

At $t=0$ we give initial data
\begin{equation}
u(x,0)=f_1(x),\quad u_t(x,0)=f_2(x). \label{5.1.3}
\end{equation}
Also, $f_1,~f_2$ are smooth functions with compact support.

We want to show that the problem is well posed. Fourier
transform with respect to $x$ gives us
\[
\hat u_{tt} = -\hat P_0(\omega)\hat u,\quad \hat P_0(\omega)=A_1\omega_1^2 + \sum_{j=2}^r
B_j\omega_j^2,\quad \hat P_0(\omega)=\hat P_0^*(\omega)>0
\]
and
\[
\hat u(\omega,0)=\hat f_1(\omega), \quad \hat u_t(\omega,0)=\hat f_2(\omega),
\]
Introducing a new variable
\[
\hat u_t = \hat P_0^{1/2} \hat v,
\]
gives us
\[
\begin{pmatrix}\hat u\cr \hat v\end{pmatrix}_t = \begin{pmatrix}0 & \hat
P_0^{1/2}\cr -\hat P^{1/2}_0 & 0\end{pmatrix}
\begin{pmatrix}\hat u\cr\hat v\end{pmatrix}.
\]
Therefore we obtain
\[
\frac{\partial}{\partial t}\Bigl(|\hat u|^2 + |\hat v|^2\Bigr) = 0,
\]
i.e.
\[
\|\hat u(\cdot,t)\|^2 + \|\hat v(\cdot,t)\|^2 = \|\hat u(\cdot,0)\|^2 + \|\hat
v(\cdot,0)\|^2.
\]
This energy estimate shows that the Cauchy problem is well posed.

%%%%%%%%%%%%%%%%%%%%%%%%%%%%%%%%%%%%%%%%%%%%%%%%%%%%%%%%%%%%%%%%%%%%%%%%
% old sec_5-2.tex
%%%%%%%%%%%%%%%%%%%%%%%%%%%%%%%%%%%%%%%%%%%%%%%%%%%%%%%%%%%%%%%%%%%%%%%%
\subsection{The resolvent equation}\label{sec_5-2}
Consider the Cauchy problem for the inhomogeneous system (\ref{5.1.1}) with $F(x,t)\in
\xCinfty_0$ and with homogeneous initial data $f_1=f_2=0.$
\begin{equation}
u_{tt}=P_0(D) u+F. \label{5.2.1}
\end{equation}
Fourier transform with respect to $x$ and Laplace
transform with respect to time gives us the resolvent equation
\begin{equation}
\Bigl(s^2I+|\omega|^2\hat P_0(\omega')\Bigr)\hat u=\hat F,\quad s=i\xi+\eta,~
\eta>0,~\omega'=\omega/|\omega|.  \label{5.2.2}
\end{equation}
Since $P_0=P_0^*>0,$ there is a unitary transformation which transforms (\ref{5.2.2})
to diagonal form
\[
\Bigl(s^2+|\omega|^2\mu_j\Bigr)\tilde u_j=\tilde F_j,~j=1,2,\ldots, n,\quad
\hat u=S_1\tilde u,~\hat F=S_1\tilde F,
\]
i.e.,
\[
(s+i|\omega|\sqrt{\mu_j}) (s-i|\omega|\sqrt{\mu_j})\tilde u_j=\tilde F_j.
\]
Without restriction we can assume that $\xi>0.$ Then
\[\begin{split}
|s+i|\omega|\sqrt{\mu_j}|&=\sqrt{|\xi+|\omega|\sqrt{\mu_j}|^2+\eta^2}\ge
\sqrt{|s|^2+|\omega|^2\mu_j},\\
|s-i|\omega|\sqrt{\mu_j}|&=\sqrt{|\xi-|\omega|\sqrt{\mu_j}|^2+\eta^2}\ge
\eta =\re s.
\end{split}
\]
Therefore
\[
|\tilde u_j|\le \frac{|\tilde F_j|}{ \sqrt{|s|^2+|\omega|^2\mu_j}\,\re s}.
\]
Choosing Im$\,s=i|\omega|\sqrt{\mu_j}$ shows that the estimate is sharp.

We have proved
\begin{thrm}\label{theorem5.2.1}
There is a constant
$K$ which does not depend on $\omega$ such that the resolvent estimate
\begin{equation}
|\hat u(\omega,s)|\le \frac{K|\hat F|}{ \sqrt{|s|^2+|\omega|^2}\, \re s}
\label{5.2.3}
\end{equation} 
holds.
\end{thrm}

The last estimate shows that we `gain' one derivative, i.e., if the forcing $\in
\xHn{p},$ then the solution $\in \xHn{{p+1}}.$ Therefore we can prove that the
Cauchy problem is stable against lower order perturbations.
\begin{thrm}\label{theorem5.2.2}
Consider the Cauchy problem with homogeneous initial data for
\[
u_{tt}=\left(P_0(D)+P_1(D)\right)u+F.
\]
Here $P_1(D)$ represents a general first order operator. There is an $\eta_0>0$ such that
 the estimate (\ref{5.2.3}) holds for $\eta=\re s >\eta_0.$
\end{thrm}
\begin{proof} We consider $P_1(D)u$ as part of the forcing. Then (\ref{5.2.3}) gives us
\[
|\hat u(\omega,s)|\le \frac{K}{ \re s} \frac{|P_1(i\omega,s)\hat u|}{\sqrt{|s|^2+|\omega|^2}}
+ \frac{K}{\re s} \frac{|\hat F|}{\sqrt{|s|^2+|\omega|^2}} .
\]
Since $\frac{|P_1(i\omega,s)|}{\sqrt{|s|^2+|\omega|^2}}$ is uniformly bounded,
we choose $\eta_0$ such that
\[ 
\frac{K}{ \eta_0}\frac{|P_1(i\omega,s)\hat u|}{\sqrt{|s|^2+|\omega|^2}} \le
\frac{1}{2} |\hat u(\omega,s)|.
\]
Then the desired estimate follows.
\end{proof}

We can write the resolvent equation also as a first order system.  We Fourier
transform  (\ref{5.2.1}) with respect to $x_-=(x_2,\ldots,x_r)$ and Laplace transform
it with respect to $t.$ Then we obtain
\begin{gather}
A_1D_1^2 \hat u=\left(s^2 I+B(\omega_-)\right)\hat u-\hat F,\quad
B(\omega_-)=\sum_{j=2}^r B_j\omega_j^2, \label{5.2.4}\\
\hat u=\hat u(x_1,\omega_-,s) \in \xLtwo(0\le x_1< \infty).\nonumber
\end{gather}
Since $A_1>0,$ there is a constant $\sigma >0$ such
that $A_1=A_1^*\ge \sigma I>0.$ Introducing a
new variable by
\[ D_1\hat u=\sqrt{|s|^2+|\omega_-|^2}\, \hat v,  \]
we obtain the first order system
\begin{equation}\label{5.2.5}
D_1\begin{pmatrix}\hat u\cr \hat v\end{pmatrix} =M(s,\omega_-)
\begin{pmatrix}\hat u\cr \hat v\end{pmatrix}
+ \frac{1}{\sqrt{|s|^2+|\omega_-|^2}} \begin{pmatrix}0\cr -A_1^{-1}\hat
F\end{pmatrix},
\end{equation}
where
\[
M=M(s,\omega_-)=\begin{pmatrix} 0 & \sqrt{|s|^2+|\omega_-|^2}I \cr
\frac{A_1^{-1}\left(s^2I+B(\omega_-)\right)}{\sqrt{|s|^2+|\omega_-|^2}} & 0 \end{pmatrix}.
\]
The eigenvalues $\kappa$ of $M$ are solutions of
\begin{equation}
A_1\kappa^2\varphi_0=\left(B(\omega_-)+s^2 I\right)\varphi_0. \label{5.2.6}
\end{equation}
\begin{lmm}\label{lemma5.2.1}
For $\re s>0,$ there are no $\kappa_j$ with
$\re\kappa_j=0.$ Also, there are exactly $n$ eigenvalues, counted according to their
multiplicity with $\re\kappa <0$ and, therefore, $n$ eigenvalues with
$\re\kappa>0.$
\end{lmm}
\begin{proof} Assume there exists a $\kappa=i\omega_1$ which is purely imaginary.
Then, by (\ref{5.2.2}),
\begin{equation}
\left(s^2I-\hat P_0(i\omega_1,i \omega_-)\right)\varphi_0=0 \label{5.2.7}
\end{equation}
has a nontrivial solution, i.e., $s^2$ is an eigenvalue of $\hat P_0(i\omega).$
$\hat P_0(i\omega)<0$ implies that $s^2$ is real and negative which is a
contradiction with $\re s>0.$

The solutions of (\ref{5.2.6}) are continuous functions of $\omega_-.$ Therefore the
number of $\kappa$ with $\re\kappa<0$ does not depend on $\omega_-$ and we can
assume that $\omega_-=0.$ Then (\ref{5.2.6}) reduces to
\begin{equation}
(s^2 I-A_1\kappa^2)\varphi_0=0. \label{5.2.8}
\end{equation}
Since, by assumption, $A_1$ has positive eigenvalues $\mu_j$ and a 
complete system of eigenvectors, we can transform (\ref{5.2.8}) into $n$ scalar
equations
\[ (s^2-\mu_j\kappa^2)u_j=0,\quad j=1,2,\ldots,n,\]
i.e.,
\[ \kappa=\pm s/\sqrt{\mu_j}, \quad \mu_j>0,~j=1,2,\dots n,~\re s>0. \]
This proves the lemma.
\end{proof}

By Schur's lemma, there exists a unitary transformation $U=U(s,\omega_-) $ such that
\begin{equation}\label{eq:108}
U^*(s,\omega_-)M(s,\omega_-)U(s,\omega_-)=\begin{pmatrix}M_{11} & M_{12}\cr 0 &
M_{22}\end{pmatrix},
\end{equation}
where the eigenvalues $\kappa_{j1},\kappa_{j2}$ of $M_{11}$ and $M_{22}$ satisfy
$\re\kappa_{j1}<0,\re\kappa_{j2}>0,$ respectively, for $\re s >0.$ Clearly, the transformed 
equation (\ref{5.2.5}) can be solved uniquely for $\re s>0.$ 

Using (\ref{5.2.3}), we shall now derive estimates for the solutions of (\ref{5.2.5}).
To accomplish this we consider a more general forcing. We replace
\[
\frac{1}{\sqrt{|s|^2+|\omega_-|^2}} \begin{pmatrix}0\cr A_1^{-1}\hat F\end{pmatrix}\quad {\rm by}\quad
\begin{pmatrix}\hat F_1\cr \hat F_2 \end{pmatrix},\quad \hat F_j=\hat F_j(x_1,\omega_-,s).
\]
We Fourier transform (\ref{5.2.5}) with respect to $x_1$ and consider
\[\begin{split}
i\omega_1\hat u &=
\sqrt{|s|^2+|\omega_-|^2}
\hat v-\hat F_1,\\
i\omega_1 \hat v &=
\frac{A_1^{-1}\left(s^2I+B(\omega_-)\right)}{\sqrt{|s|^2+|\omega_-|^2}} \hat
u-\hat F_2.
\end{split}
\]
Eliminating $\hat v$ gives us
\[
\left(s^2I+|\omega|^2\hat P_0(\omega')\right) \hat u= i\omega_1 A_1\hat F_1+
\sqrt{|s|^2+|\omega_-|^2} A_1\hat F_2.
\]
Therefore, by (\ref{5.2.3}), we obtain the estimate
\[\begin{split}
|\hat u(\omega,s)|&\le  K \frac{|i\omega_1\hat F_1+
\sqrt{|s|^2+|\omega_-|^2} \hat F_2|}{\sqrt{|s|^2+|\omega|^2} \re s} \\
&\le K\frac{|\hat F_1|+|\hat F_2|}{ \re s} , \quad \hat F_j=\hat
F_j(\omega_1,\omega_-,s).
\end{split}
\]
Eliminating $\hat u,$ we obtain the same estimate for $\hat v.$  Therefore we have proved 

\begin{lmm}\label{lemma5.2.2}
There exists a constant $K>0$ such that, for all $\omega_1,
\omega_-, s,$
\begin{equation}
\Bigl|\bigl(M(s,\omega_-)-i\omega_1I\bigr)^{-1}\Bigr|\le \frac{2K}{\re s}. \label{5.2.9}
\end{equation}
In particular, the eigenvalues $\kappa$ of $M(s,\omega_-)$ satisfy
\begin{equation}
|\re\kappa|\ge \frac{\re s}{ 2K}. \label{5.2.10}
\end{equation}
\end{lmm}

Using scaled variables
\begin{equation}
s'=s/\sqrt{|s|^2+|\omega_-|^2},\quad
\omega'=\omega/\sqrt{|s|^2+|\omega_-|^2},\quad
M'(s',\omega')=\begin{pmatrix}0 & I\cr A_1^{-1}\bigl(s^{\prime 2}I+\beta(\omega'_-)\bigr) &
0\end{pmatrix},
\label{5.2.11}
\end{equation}
we can write (\ref{5.2.9}),(\ref{5.2.10}) in the form
\begin{align}
|\left(M'(s',\omega'_-)-i\omega'_1 I\right)^{-1}|&\le \frac{2K}{ \re s'},\label{5.2.12}\\
|\re\kappa'| &\ge \frac{\re s'}{ 2K}.\label{5.2.13}
\end{align}

%%%%%%%%%%%%%%%%%%%%%%%%%%%%%%%%%%%%%%%%%%%%%%%%%%%%%%%%%%%%%%%%%%%%%%%%
% old sec_5-3.tex
%%%%%%%%%%%%%%%%%%%%%%%%%%%%%%%%%%%%%%%%%%%%%%%%%%%%%%%%%%%%%%%%%%%%%%%%
\subsection{Reduction to a first order system of pseudo-differential equations}\label{sec_5-3}
Now we consider the general halfplane problem (\ref{2.1.1})--(\ref{2.1.3})
with homogeneous initial data $f_1(x)=f_2(x)=0,$ coupled to the boundary
conditions (\ref{2.1.3}). We Fourier transform the problem with respect to
$x_-=(x_2, x_3, \dots x_r)$ and Laplace transform it with respect to $t$ and
obtain (\ref{5.2.4}) coupled to the boundary condition
\begin{equation}\label{eq:113}
\hat u_x + (C_0 s + \sum_{j=2}^r C_j\omega_j)\hat u=\hat g, \quad \hat u=\hat
u(0,\omega_-,s), \quad \hat g=\hat g(\omega_-, s)
\end{equation}

As in section \ref{sec_3} we have to assume that there are no simple wave solutions  for $\re
s>0$~ i.e., that the eigenvalue problem consisting of the homogeneous equations
(\ref{5.2.4}) and (\ref{eq:113}) have no eigenvalues $s$ with $\re s>0,$
otherwise the problem is not well posed.

Now we introduce new variables by
\begin{equation}\label{eq:114}
\hat u= A_1^{-1/2} \tilde u_1, \quad D_1\tilde u=\sqrt{|s|^2+|\omega_-|^2} \tilde
v
\end{equation}
and obtain a normalised version of (\ref{5.2.5})
\begin{equation}
D_1\begin{pmatrix} \tilde u\cr \tilde v\end{pmatrix}=\sqrt{|s|^2+|\omega_-|^2}\, \tilde M
\begin{pmatrix} \tilde u\cr \tilde v\end{pmatrix}+\frac{1}{ \sqrt{|s|^2+|\omega_-|^2}}
\begin{pmatrix}0\cr -A_1^{-1/2}\hat F \end{pmatrix} \label{5.3.3}
\end{equation}
where
\[
H(s',\omega')=A_1^{-1/2}\left(s^{\prime 2}I+B(\omega'_-)\right)A^{-1/2},\quad \tilde
M=\begin{pmatrix}0 &I\cr H & 0\end{pmatrix},
\]
and 
$ s'=s/ \sqrt{|s|^2+|\omega_-|^2},~ \omega'_-= \omega_-/\sqrt{|s|^2+|\omega_-|^2} $ are
scaled variables.

Instead of (\ref{5.2.6}) we obtain now
\begin{equation}\label{5.3.4}
\begin{split}
\psi&= \kappa' \varphi,\\
H(s',\omega')\varphi &=\kappa'\psi,
\end{split}
\end{equation}
i.e.,
\begin{equation}
\kappa^{\prime 2}\varphi=A_1^{-1/2}\left(s^{\prime 2}I+B(\omega'_-)\right)
A_1^{-1/2}\varphi, \quad \kappa'=\frac{\kappa}{\sqrt{|s|^2 + |\omega_-|^2}}
\label{5.3.5}
\end{equation}
(\ref{5.3.5}) is a normalized form of (\ref{5.2.6}).

Lemmas \ref{lemma5.2.1} and (\ref{lemma5.2.2}) are crucial because they
garantee that we can use the classical theory (see lemmas 2.1 and 2.2 in [1]). In
particular, as in section \ref{sec_3}, away from any eigenvalue or generalised
eigenvalue $s,$ with $\re s>0,$ the strong estimate of definition
\ref{definition2.1.3} holds. Thus we need only to consider the estimates in the
neigbourhood of generalised eigenvalues.

Also, for $\re s>0$ we can use Schur's lemma to transform $\tilde M$ to the upper triangular form
(\ref{eq:108}) separating the eigenvalues with $\re \kappa_j<0$ and $\re \kappa_j>0$
respectively. Then we use the technique of section \ref{sec_3} to estimate the solutions. This
becomes particularly simple if there is a standard energy estimate and we need only to show that the problem
is {\it Boundary Stable}.

Finally, we want to show that in the neighbourhood of generalised eigenvalues
our problem behaves like wave equations.

We shall now derive a normal form of $\tilde M$ for $s_0'=i\xi_0'.$  We have
\begin{lmm}\label{lemma5.3.3}
Let $s_0'=i\xi_0',$ $\omega_{-0}'$ be a generalised eigenvalue.
$H(\omega_0',i\xi_0')$ is symmetric and there is a unitary transformation $U$
such that
\begin{equation}
U^*HU=\begin{pmatrix}H_{11} & 0\cr 0 & H_{22}\end{pmatrix},\label{5.3.8}
\end{equation}
where
\begin{align}
H_{11}&=\begin{pmatrix} \kappa'^2_1&&& 0\cr & \kappa'^2_2&&\cr &&\ddots &\cr 0
&&&\kappa'^2_m\end{pmatrix}
,\quad \kappa'^2_1\ge \kappa'^2_2\ge \ldots \ge \kappa'^2_m >0, \label{5.3.9}\\
H_{22}&=\begin{pmatrix} \kappa'^2_{m+1}&&& 0\cr & \kappa'^2_{m+2}&&\cr
 &&\ddots &\cr 0 &&&\kappa'^2_n\end{pmatrix}
,\quad 0\ge\kappa'^2_{m+1}\ge \kappa'^2_{m+2}\ge \ldots \ge \kappa'^2_n.
\label{5.3.10}
\end{align}
where $\kappa'_j=\pm \kappa'_j(\omega'_{-0},i\xi'_0)$ are the eigenvalues of
(\ref{5.3.5}).
 $\kappa'^2_{m+1}=0$ if and only if 
\[
\xi_0'^2=\beta_j|\omega'_{-0}|^2\quad \hbox{is an eigenvalue of}~B(\omega'_{-0}),
\quad 0<\beta_{min}\le \beta_j\le \beta_{max},~j=1,2,\ldots,n.
\]
Also,
\begin{equation}
{\rm if}\quad \xi_0'^2 <\beta_{min} |\omega'_{-0}|^2,\quad \hbox{ then all}~
\kappa'^2_j >0. \label{5.3.11}
\end{equation}
\end{lmm}
\begin{proof} Since $H(\omega'_{-0},i\xi'_0)$ is symmetric, we obtain
(\ref{5.3.8})--(\ref{5.3.10}).

If $\xi'^2_0>\beta_{max}|\omega'_{-0}|^2,$ then $H$ is negative definite and all $\kappa_j'^2<0.$
Correspondingly, if $\xi'^2_0<\beta_{min}|\omega'_{-0}|^2,$ then all $\kappa'^2_j >0.$
If $\kappa'^2=0,$ then (\ref{5.3.5}) becomes
\[
\left(B(\omega'_{-0})-\xi'^2_0I\right)\varphi=0,
\]
i.e., $\xi'^2_0$ must be an eigenvalue of $B(\omega'_{-0}).$ Clearly, the reverse is also true. If
$\xi'^2_0$ is an eigenvalue of $B(\omega'_{-0}),$ then there is a $\kappa'^2=0.$ This proves the lemma.
\end{proof}
To simplify the arguments we shall make a strong assumption which we shall relax at the end
of the section.
\begin{ssmptn}\label{assumption5.3.2} The eigenvalues $\kappa'^2_j$ are distinct.
\end{ssmptn}

In this case we can choose $U=U(\omega'_-,\xi')$ as a smooth function of
$\omega'_-,\xi'$ in a neighborhood of $\omega'_{-0},\xi'_0.$ Also there is a
constant $d_0>0$ such that, in the whole neighborhood,
\begin{equation}
|\kappa'^2_j-\kappa'^2_i|\ge d_0\quad \hbox{for all $i,j$ with}\quad i\ne j.
\label{5.3.12}
\end{equation}

Finally, we make the perturbation $s'=i\xi'+\eta',~-\eta'_0\le \eta'\le
\eta'_0,~\eta'_0\ll 1.$
We have
\begin{equation}\label{5.3.13}
\begin{split}
& U^*(\omega'_-,\xi')H(\omega'_-,s')U(\omega'_-,\xi')=
 U^*(\omega'_-,\xi')H(\omega'_-,i\xi')U(\omega'_-,\xi')\\
&\quad + U^*(\omega'_-,\xi')A^{-1}U(\omega'_-,\xi')(2i\xi'\eta'+\eta'^2).
\end{split}
\end{equation}
Since $U^*A^{-1} U$ is strictly positive definite, its diagonal elements $a_{jj}>0$ are positive. 
By using well known algebraic results (Gershgorin's theorem) this gives us
\begin{lmm}\label{lemma5.3.4}
For sufficiently small $\eta'_0$ which depends only on $A^{-1}$ and
(\ref{5.3.12}),
there exists a smooth nonsingular transformation $S=I+\eta'\, S_1(\omega'_-,s')$ such that
\begin{equation}\label{5.3.14}
\tilde H(\omega'_-,s')=S^{-1}U^*HUS=\begin{pmatrix}H_{11} & 0\cr 0 &
H_{22}\end{pmatrix} +
2i\xi'\eta' \begin{pmatrix} \tilde a_{11} &&0\cr &\ddots &\cr 0 && \tilde
a_{nn}\end{pmatrix},
\end{equation} 
where
$\tilde a_{jj}=a_{jj}+{\cal O}(\eta')>0$ and $H_{11},H_{22}$ are given as before
but with distinct eigenvalues.
\end{lmm}

Now we can construct the normal form for the resolvent equation (\ref{5.3.3}).
We introduce new variables by
 \[
\tilde u=A_1^{-1/2} US\,\tilde{\tilde{u}},\quad \tilde v=A_1^{-1/2}
US\,\tilde{\tilde{v}}
\]
and after a permutation we obtain
\begin{thrm}\label{theorem5.3.1}
In a neighborhood of $\omega'_{-0},\xi'_0$ the resolvent equation
(\ref{5.3.3}) can be transformed smoothly into
\begin{equation}
D_1
\begin{pmatrix}\tilde{\tilde{u}}\cr \tilde{\tilde{v}}\end{pmatrix}=
\sqrt{|s|^2+|\omega_-|^2}\, \begin{pmatrix}0 & I \cr \tilde H(\omega'_-,s') &
0\end{pmatrix}
\begin{pmatrix}\tilde{\tilde{u}}\cr \tilde{\tilde{v}}\end{pmatrix}+
\frac{1}{\sqrt{|s|^2+|\omega_-|^2}}\begin{pmatrix}0\cr \tilde{\tilde{F}}\end{pmatrix}
\label{5.3.15}
\end{equation} 
where $\tilde{\tilde{F}}= -S^{-1}U^*A^{1/2}_1\hat F.$ By (\ref{5.3.14}), the system
(\ref{5.3.15}) is
composed of $2\times 2$ systems
\begin{equation}
\begin{pmatrix}\tilde{\tilde{u}}_j\cr \tilde{\tilde{v}}_j\end{pmatrix}_x=
\sqrt{|s|^2+|\omega_-|^2}\, \begin{pmatrix}0 & 1 \cr \tilde \kappa'^2_j(\omega',s') &
0\end{pmatrix}
\begin{pmatrix}\tilde{\tilde{u}}_j\cr \tilde{\tilde{v}}_j\end{pmatrix}+
\frac{1}{\sqrt{|s|^2+|\omega_-|^2}}\begin{pmatrix}0\cr \tilde{\tilde{F}}_j\end{pmatrix}
\label{5.3.16}
\end{equation} 
where
\[
\tilde \kappa'^2_j(\omega,s)=\kappa'^2_j(\omega'_-,i\xi')+2i\tilde a_{jj}\xi'\eta'.
\]
\end{thrm}

Clearly, the $2\times 2$ blocks have the form (\ref{3.2.4}) of section
\ref{sec_3}.

There are no difficulties to generalize the results to the case that the
$\kappa_j^2$ have constant multiplicity (this was done in [5] for first orders systems).

%%%%%%%%%%%%%%%%%%%%%%%%%%%%%%%%%%%%%%%%%%%%%%%%%%%%%%%%%%%%%%%%%%%%%%%%
% old sec_4-2.tex (Anders' file)
%%%%%%%%%%%%%%%%%%%%%%%%%%%%%%%%%%%%%%%%%%%%%%%%%%%%%%%%%%%%%%%%%%%%%%%%
\section{Numerical experiments}\label{sec_5}

In this section we numerically solve the strip problem for the scalar wave equation
\begin{equation}
u_{tt}=u_{xx}+u_{yy}+F(x,y,t),\quad  0 \leq x \leq 1,\ 0 \leq y \leq 1,\ t\ge
0,\label{eq:scalar-wave-eqn}
\end{equation}
with 1-periodic solutions in the $y$-direction,
\begin{equation}
u(x,y,t)=u(x,y+1,t),\quad 0\leq x\leq 1,\ t\geq 0,\label{eq:periodic-cond}
\end{equation}
subject to initial conditions,
\begin{equation}
u(x,y,0) = f_1(x,y),\quad u_t(x,y,0)=f_2(x,y),\quad  0 \leq x \leq 1,\ 0 \leq y
\leq 1,\label{eq:initial-cond}
\end{equation}
and boundary conditions
\begin{alignat}{2}
u_x - b\, u_y &= g_0(y,t),\quad x=0,\ 0\leq y \leq 1,\ t\ge 0,\label{eq:bc1}\\
u_x &= g_1(y,t),\quad x=1,\ 0\leq y \leq 1,\ t\ge 0.\label{eq:bc2}
\end{alignat}
Here $b$ is a constant.  We are interested in the three cases $b=0$, $b$ real, and $b$
purely imaginary, i.e.,
\[
b=i\beta,\quad\mbox{$\beta$ real.}
\]
To solve the latter problem using real arithmetic, we introduce real-valued functions
$u^{(1)}$ and $u^{(2)}$ such that
\begin{equation}\label{eq:real}
u = u^{(1)} + i u^{(2)}.
\end{equation}
Inserting (\ref{eq:real}) into  (\ref{eq:scalar-wave-eqn}) leads to the system of scalar wave equations,
\begin{alignat}{2}
u^{(1)}_{tt}&=u^{(1)}_{xx}+u^{(1)}_{yy}+{\rm Re}\,F(x,y,t),\quad  0 \leq x \leq 1,\quad 0
\leq y \leq 1,\quad t\ge 0,\label{eq:scalar-wave-eqna}\\ 
u^{(2)}_{tt}&=u^{(2)}_{xx}+u^{(2)}_{yy}+{\rm Im}\,F(x,y,t),\quad  0 \leq x \leq 1,\quad 0
\leq y \leq 1,\quad t\ge 0.\label{eq:scalar-wave-eqnb}
\end{alignat}
Boundary condition (\ref{eq:bc1}) can be written as
\begin{alignat}{2}
u^{(1)}_x + \beta \, u^{(2)}_y &= {\rm Re}\,g_0(y,t),\quad x=0,\ 0\leq y\leq 1,\ t\ge 0,\label{eq:3re}\\
u^{(2)}_x - \beta \, u^{(1)}_y &= {\rm Im}\,g_0(y,t),\quad x=0,\ 0\leq y\leq 1,\ t\ge 0,\label{eq:3im}
\end{alignat}
which is of the form (\ref{4.1.3}) with $b_1=\beta$ and $b_2=-\beta$.

We introduce a grid with grid size $h=1/(N-1)$,
\[
x_j=(j-1)h,\quad j=0,1,\ldots,N+1,\quad y_k = k h,\quad k=0,1,\ldots,N.
\]
Time is discretized on a uniform grid with time step
$\delta_t>0$, $t_n=n \delta_t$, $n=0,1,2,\ldots$ and we denote a grid function by
\[
v^n_{j,k} = v(x_j, y_k, t_n).
\]
The standard divided difference operators are defined by
\[
D_{+x} v^n_{j,k} = \frac{v^n_{j+1,k} - v^n_{j,k}}{h},\quad D_{-x} v^n_{j,k} = D_{+x}
v^n_{j-1,k},\quad D_{0x} = \frac{1}{2}\left( D_{+x} + D_{-x}\right),
\] 
with corresponding notations in the $y$- and $t$-directions.

Consider the difference approximation
\begin{equation}\label{eq:da}
D_{+t}D_{-t} v^{n}_{j,k} = \left( D_{+x}D_{-x}+D_{+y}D_{-y} \right) v^{n}_{j,k} + F(x_j, y_k, t_n),
\end{equation}
subject to boundary conditions
\begin{alignat}{3}
v^{n}_{j,0} - v^{n}_{j,N-1} &= 0,\quad &j=0,1,\ldots,N+1,\\
v^{n}_{j,N} - v^{n}_{j,1} &= 0,\quad &j=0,1,\ldots,N+1,\\
D_{0x} v^{n}_{1,k} - b\, D_{0y} v^{n}_{1,k} &= g_0(y_k,t_n),\quad &k=1,2,\ldots,N-1,\label{eq:disc-bc1}\\
D_{0x} v^{n}_{N,k} &= g_1(y_k,t_n),\quad &k=1,2,\ldots,N-1,\label{eq:disc-bc2}
\end{alignat}
for $n=-1,0,1,\ldots$, and initial conditions,
\begin{equation}
v^{0}_{j,k} = f_1(x_j,y_k),\quad v^{-1}_{j,k} = f'_2(x_j,y_k),\quad j=0,1,\ldots,N+1,\ k=0,1,\ldots,N.\label{eq:d-ic}
\end{equation}
When $b=0$, the difference approximation (\ref{eq:da})-(\ref{eq:d-ic}) satisfies a discrete energy
estimate, under the Courant time step restriction
\[
\delta_t \leq C h,
\]
and is therefore stable. The energy method can not be used to show stability of the difference
approximation when $b$ is non-zero and real, or purely imaginary. However, as we shall see below,
our practical experience indicates that the approximation is stable also when $b$ is purely
imaginary. When $b$ is non-zero and real, the continuous problem is {\it Unstable}. In this case,
the difference approximation is convergent for short times.

To test the accuracy of the numerical solution, we choose the forcing functions $F$,
$g_k$, and $f_k$ such that the exact solution becomes the traveling wave
\[
u_w(x,y,t) = \sin(2\pi (x-t)) \sin(2\pi y).
\]
For the case $b=i\, \beta$, we use the exact solution
\begin{alignat*}{2}
u_w^{(1)}(x,y,t) &= \sin(2\pi (x-t)) \sin(2\pi y),\\
u_w^{(2)}(x,y,t) &= \cos(2\pi (x-t)) \cos(2\pi y).
\end{alignat*}
Table~\ref{tab:accuracy} shows the max norm of the error $u_w-v$ at different times for grid
sizes $h=10^{-2},\ 5\times 10^{-3},\ 2.5\times 10^{-3}$. All calculations used $\delta_t =
0.5 h$. The error decreases as ${\cal O}(h^2)$ for all three values of
$b$, both at time $t=1$ and $t=10$.
\begin{table}[ht]
\begin{center}
\begin{tabular}{c|c|c|c}
Case    & $h$ & $\| u_w-v \|_\infty(t=1)$ & $\| u_w-v \|_\infty(t=10)$  \\ \hline
$b=0$   & $1\times 10^{-2}$ & $7.09\times 10^{-4}$ & $4.69\times 10^{-4}$  \\
        & $5\times 10^{-3}$ & $1.76\times 10^{-4}$ & $1.18\times 10^{-4}$  \\
        & $2.5\times 10^{-3}$ & $4.42\times 10^{-5}$ & $2.97\times 10^{-5}$  \\ \hline
$b=0.5$ & $1\times 10^{-2}$   & $7.09\times 10^{-4}$ & $2.95\times 10^{-2}$  \\
        & $5\times 10^{-3}$   & $1.76\times 10^{-4}$ & $7.37\times 10^{-3}$  \\
        & $2.5\times 10^{-3}$ & $4.42\times 10^{-5}$ & $1.84\times 10^{-3}$  \\ \hline
$b=i\,0.5$ & $1\times 10^{-2}$   & $1.44\times 10^{-3}$ & $1.61\times 10^{-3}$ \\
           & $5\times 10^{-3}$   & $3.61\times 10^{-4}$ & $4.02\times 10^{-4}$ \\
           & $2.5\times 10^{-3}$ & $9.03\times 10^{-5}$ & $1.00\times 10^{-4}$
\end{tabular}
\end{center}
\caption{Max error in the solution as function of the grid size with a
  traveling wave exact solution.}\label{tab:accuracy}
\end{table}

\begin{figure}[ht]
  \begin{center}
    \includegraphics[width=0.9\textwidth]{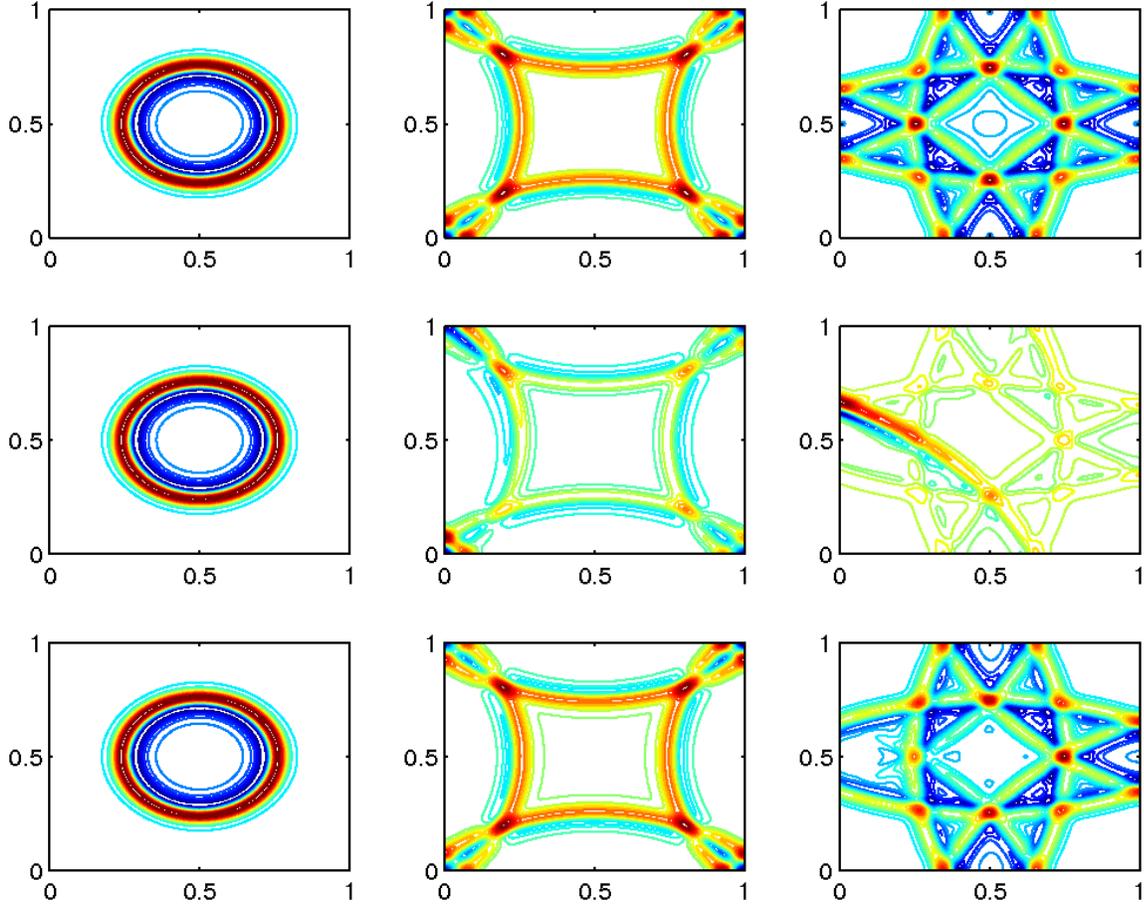}
  \end{center}
  \caption{The solution at times $t=0.25$ (left column), $t=0.75$ (middle column), and
  $t=1.25$ (right column) for $b=0$ (top row), $b=0.5$ (middle row), and $b=i\, 0.5$
  (bottom row). The bottom row is showing the real part of the solution ($u^{(1)}$).}
  \label{fig:cases}
\end{figure}

\begin{figure}[ht]
  \begin{center}
    \includegraphics[width=0.9\textwidth]{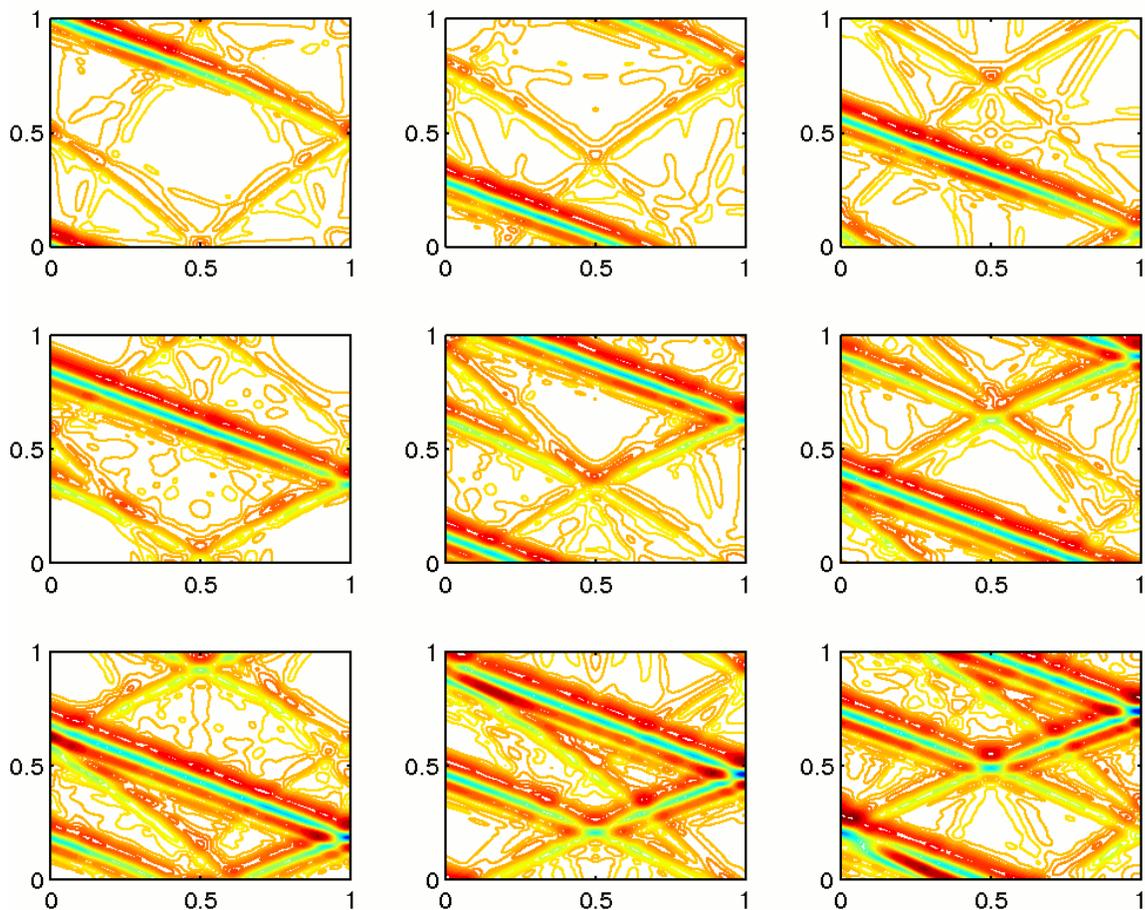}
  \end{center}
  \caption{The solution of the unstable case ($b=0.5$) at times $t=2.5$ to
  $t=4.5$ in increments of $0.25$, starting in the top left sub-figure and progressing
  row-wise to the bottom right sub-figure, e.g. $t=2.75$ is in the middle column of the
  top row. }
  \label{fig:evol}
\end{figure}

To illustrate how the $b$-coefficient in the boundary condition influences the solution, we study the
evolution from an initial Gaussian pulse,
\[
f_1(x,y) = f'_2(x,y) = e^{-(x-0.5)^2/L^2 - (y-0.5)^2/L^2},\quad L=0.03,
\]
with homogeneous interior and boundary forcing functions, $F=0$, $g_0=0$, $g_1=0$. In
these calculations we use the grid size $h=5\times 10^{-3}$ and time step $\delta_t =
2.5\times 10^{-3}$. The evolution between times $t=0.25$ and $t=1.25$ for the three cases
$b=0$, $b=0.5$, and $b=i\,0.5$ is shown in Figure~\ref{fig:cases}. The solution initially
propagates outwards towards the boundary (first column). Before the pulse reaches the
boundary, the solutions are identical for all three cases. The middle column in
Figure~\ref{fig:cases} shows the solutions after the initial pulse has reached the
boundary, at $t=0.75$. At this time, there are only minor differences between the three
solutions. The influence of the boundary condition is becoming more obvious in the right
column, corresponding to $t=1.25$. While the differences between $b=0$ and $b=i\,0.5$ are
still small and located near the left boundary, the case $b=0.5$ has developed a structure near
that boundary which is not present for $b=0$, or $b=i\, 0.5$.
For later times, the boundary structure develops into a diagonal streak which extends
further and further into the domain until it gets reflected by the Neumann condition on
the right boundary, see Figure~\ref{fig:evol}. At later times, the reflected streak
develops a new streak which grows into the domain in the same way. The solution is
eventually dominated by these streaks which appear to propagate in the direction
$y+0.5x$. Note that the propagation direction $y+0.5 x$ is consistent with eigenfunction
(\ref{eq:eigenfcn4}) since $b=0.5$, see case 4) in section \ref{sec_3-1}.

\begin{figure}[ht]
  \begin{center}
    \includegraphics[width=0.9\textwidth]{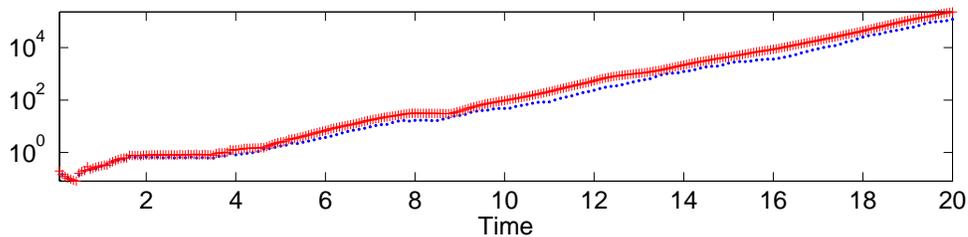}
  \end{center}
  \caption{The max norm of the solution for $0\leq t\leq 20$ for the case
$b=0.5$, starting from a Gaussian pulse. The blue dots correspond to grid size
$h=10^{-2}$ and the red crosses have $h=5\times 10^{-3}$. }
  \label{fig:long-time}
\end{figure}

It is also interesting to monitor the max norm of the solution for longer times when
$b=0.5$, see Figure~\ref{fig:long-time}. Note that the solution grows exponentially with
time, illustrating the {\it Unstable} nature of this boundary condition. Also note that the solution
is slightly larger on the finer grid. This behavior agrees with the predicted exponential
growth proportional to $|\omega|^t$, because higher values of $|\omega|$ are captured on
the finer grid. Note, however, that this growth is not due to numerical instabilities
because the accuracy test shows second order convergence, at least up to $t=10$, see
Table~\ref{tab:accuracy}.

\begin{figure}[htb]
  \begin{center}
    \includegraphics[width=0.7\textwidth]{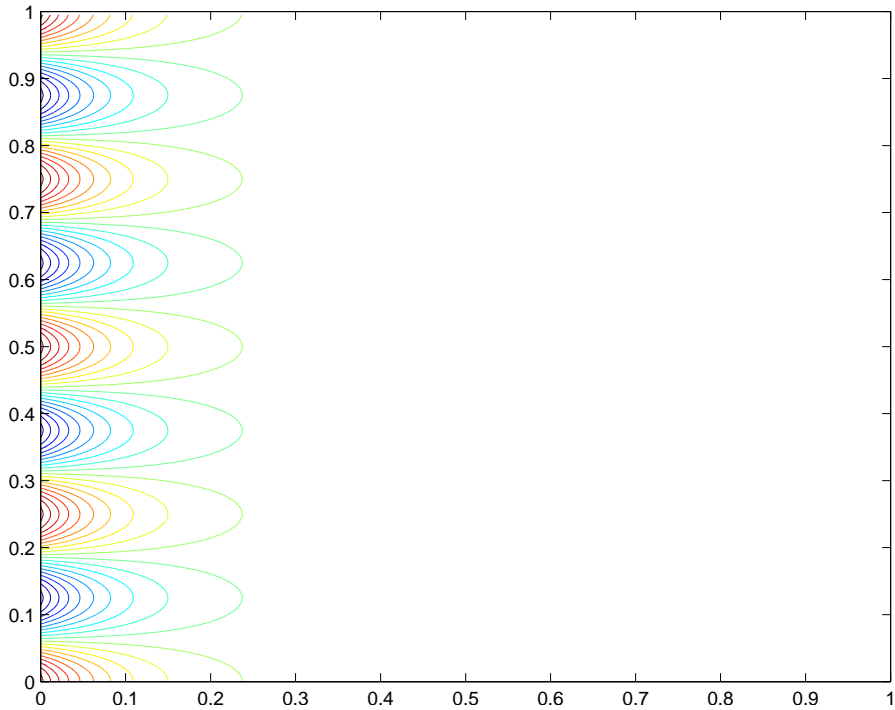}
  \end{center}
  \caption{The real part of the initial data for the surface wave with
$\beta=0.5$ and $\omega_0 = 8\pi$. }
  \label{fig:initial}
\end{figure}

To more clearly see the difference between the cases $b=0$ and $b=i\, \beta$ we take
$F=0$, $g_0=g_1=0$ and change the initial data to trigger a surface wave,
\[
f_1(x,y) = u_s(x,y,0),\quad f'_2(x,y) = u_s(x,y,-\delta_t),
\]
where
\begin{equation}\label{eq:surface-wave}
u_s(x,y,t) = e^{-|\beta\omega_0|x}\left[\cos\left(\omega_0(y - \sqrt{1-\beta^2}\,t)\right) + i\,
\sin\left(\omega_0 (y - \sqrt{1-\beta^2}\,t)\right) \right],\quad \beta\omega_0>0. 
\end{equation}
This wave decays exponentially away from the $x=0$ boundary with a harmonic oscillation in
$y$, see Figure~\ref{fig:initial}. The surface wave propagates in the positive
$y$-direction with a wave speed proportional to $\sqrt{1-\beta^2}$. As $\beta\to 0$, the
surface wave decays slower and slower in the $x$-direction. In the limit $\beta=0$, the
amplitude of the wave is constant in $x$ which corresponds to one-dimensional wave
propagation in the $y$-direction, consistent with the limiting boundary condition
$u_x=0$. There are no numerical difficulties in this limit.

\begin{figure}[ht]
  \begin{center}
    \includegraphics[width=0.9\textwidth]{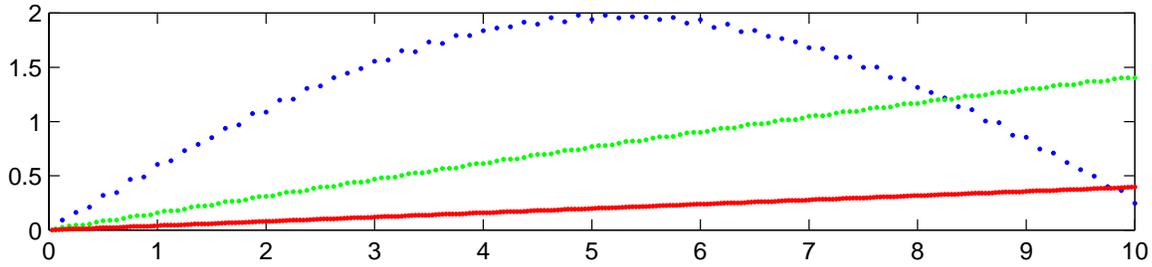}
  \end{center}
  \caption{The max norm of the error as function of time for a surface wave with
    $\beta=0.9$ computed on a grid with $h=10^{-2}$ (blue), $h=5\times 10^{-3}$ (green),
    and $h=2.5\times 10^{-3}$ (red).}
  \label{fig:phase}
\end{figure}

The case $|\beta|\to 1$ is more difficult to solve numerically. Here we study
$0.5 \leq \beta < 1$ and we use (\ref{eq:surface-wave}) as an approximation
of the exact solution ($u_s$ is only exponentially small at $x=1$ and
does not exactly satisfy the boundary condition at that boundary). To
make sure the amplitude of the surface wave is negligible at the $x=1$ boundary,
we choose
\[
\omega_0 = 8\pi,\quad e^{-|\beta\omega_0|} = e^{-4\pi} \approx 3.48\times 10^{-6},\quad
\beta=0.5.
\]
In Table~\ref{tab:surface} we show the max norm of the error $u_s-v$ for different values
of $\beta$. The case $\beta=0.5$ shows second order convergence, both at time $t=1$ and
$t=10$. As can be expected in wave propagation problems, the error is dominated by the
phase error, which explains why it is about 10 times larger at $t=10$ compared to
$t=1$. For $\beta=0.9$, the error still converges to second order accuarcy at time $t=1$,
but shows an unexpected pattern at time $t=10$. Here the error is larger for the
intermediate grid size $h=5\times 10^{-3}$ than for the coarse grid size $h=10^{-2}$. This
behavior is explained by studying the time history of the error, see
Figure~\ref{fig:phase}. For $h=10^{-2}$, the max error occurs at $t\approx 5.5$ when the
numerical solution is about 180 degrees out of phase with the exact solution. At later
times the error in the numerical solution decreases because it is between 180 and 360
degrees out of phase. The grid with $h=5\times 10^{-3}$ is barely fine enough to capture
the solution at time $t=10$ because the phase error exceeds 90 degrees. As a result we
don't see the expected second order convergence when the grid is refined to $h=2.5\times
10^{-3}$. However, the error at $t=10$ is about 10 times larger than at $t=1$ for the
finest grid, which indicates that this resolution is adequate for $\beta=0.9$. The
situation is even more dire for $\beta=0.99$. Here the errors at time $t=1$ show a simular
behavior as at $t=10$ for $\beta=0.9$, so only the finest grid provides adequate
resolution at $t=1$. At $t=10$, the error displays a completely erratic behavior with the
largest error for the finest grid. An even finer grid would be necessary to obtain an
accurate solution at $t=10$, when $\beta=0.99$.
\begin{table}[ht]
\begin{center}
\begin{tabular}{c|c|c|c}
Case         & $h$                 & $\| u_s-v \|_\infty(t=1)$ & $\| u_s-v \|_\infty(t=10)$ \\ \hline
$\beta=0.5$  & $1\times 10^{-2}$   & $2.44\times 10^{-2}$ & $2.38\times 10^{-1}$  \\
             & $5\times 10^{-3}$   & $6.35\times 10^{-3}$ & $6.19\times 10^{-2}$  \\
             & $2.5\times 10^{-3}$ & $1.60\times 10^{-3}$ & $1.56\times 10^{-2}$  \\ \hline
$\beta=0.9$  & $1\times 10^{-2}$   & $6.04\times 10^{-1}$ & $2.46\times 10^{-1}$ \\
             & $5\times 10^{-3}$   & $1.58\times 10^{-1}$ & $1.40\times 10^{0}$  \\
             & $2.5\times 10^{-3}$ & $4.00\times 10^{-2}$ & $3.95\times 10^{-1}$  \\ \hline
$\beta=0.99$ & $1\times 10^{-2}$   & $1.67\times 10^{0}$ & $1.37\times 10^{0}$ \\
             & $5\times 10^{-3}$   & $5.20\times 10^{-1}$ & $1.81\times 10^{-1}$ \\
             & $2.5\times 10^{-3}$ & $1.44\times 10^{-1}$ & $1.47\times 10^{0}$
\end{tabular}
\end{center}
\caption{Max error in the solution as function of the grid size when the exact solution is
  the surface wave $u_s(x,y,t)$.}\label{tab:surface}
\end{table}

So why is it so hard to calculate an accurate numerical solution as $|\beta|\to 1$? The
spatial resolution in terms of grid points per wave length only depends on
$\omega_0$. With $\omega_0=8\pi$, the wave length is $1/4$ and grid sizes $h=10^{-2},\
5\times 10^{-3},\ 2.5\times 10^{-3}$ correspond to 25, 50, and 100 grid points per wave
length, respectively. The exponential decay in the $x$-direction only depends weakly on
$\beta$ and never exceeds $e^{-|\omega_0|x}$ for $\beta<1$. Hence the solution varies on
the same length scale in the $x$- and $y$-directions. Furthermore, the temporal resolution
in terms of time steps per period only improves as $|\beta|\to 1$ because the wave speed
goes to zero in this limit. We conclude that the numerical difficulties are not due to
poor resolution of the solution.

To further analyze the cause of the poor accuracy in the numerical solution for $|\beta|\to
1$, we decompose the problem (\ref{eq:scalar-wave-eqn})-(\ref{eq:bc2}) into two parts,
\[
u(x,y,t) = U(x,y,t) + u'(x,y,t),
\]
such that $U$ satisfies a doubly periodic problem on an extended domain,
\[
U_{tt}=U_{xx}+U_{yy}+\tilde{F}(x,y,t),\quad  -1 \leq x \leq 2,\ 0 \leq y \leq 1,\ t\ge 0,
\]
subject to initial conditions,
\[
U(x,y,0) = \tilde{f}_1(x,y),\quad U_t(x,y,0)=\tilde{f}_2(x,y),\quad  -1 \leq x \leq 2,\ 0 \leq y \leq 1,
\]
and periodic boundary conditions
\begin{alignat*}{2}
U(x,y,t)&=U(x,y+1,t),\quad -1\leq x\leq 2,\ t\geq 0,\\
U(x,y,t)&=U(x+3,y,t),\quad 0\leq y\leq 1,\ t\geq 0.
\end{alignat*}
The interior forcing function and the initial data can be smoothly extended to become 3-periodic in the
$x$-direction, without changing them on the original domain,
\[
\tilde{F}(x,y,t)=F(x,y,t),\quad \tilde{f}_k(x,y)=f_k(x,y),\quad 0\leq x\leq 1,\ 0 \leq y
\leq 1,\ t\geq 0.
\]
The problem for $U$ is independent of the $b$-coefficient in the boundary condition and
can easily be solved numerically.

The difference $u'=u-U$ satisfies the scalar wave equation (\ref{eq:scalar-wave-eqn})-(\ref{eq:bc2}) with homogeneous
interior forcing, homogeneous initial data, but inhomogeneuos boundary conditions,
\begin{alignat}{2}
u'_x - b\, u'_y &= g'_0(y,t)\quad x=0,\ 0\leq y \leq 1,\ t\ge 0,\label{eq:bc-trouble} \\
u'_x &= g'_1(y,t)\quad x=1,\ 0\leq y \leq 1,\ t\ge 0,
\end{alignat}
The boundary forcing functions depend on $U$ according to
\begin{alignat*}{3}
g'_0(y,t) &= g_0(y,t) - \left( U_x(0,y,t) - b\, U_y(0,y,t) \right),\quad & 0\leq y \leq
1,\ t\ge 0,\\ g'_1(y,t) &= g_1(y,t) - U_x(1,y,t),& 0\leq y \leq 1,\ t\ge 0.
\end{alignat*}
The corresponding half-plane problems were analyzed in section \ref{sec_3-2}. The accuracy problems
are unlikely to arise from the Neumann boundary condition at $x=1$ since it is independent
of the $b$-coefficient. However, the half-plane problem subject to (\ref{eq:bc-trouble})
satisfies the estimates of Theorem~\Rref{theorem3.2.1}. Here, $b=i\,\beta$
corresponds to case 2), and
estimate (\ref{3.2.16}) shows that the Laplace-Fourier transform of $u'$ satisfies
\begin{equation}\label{eq:gen-eig-est}
|\tilde{u}'(0,\omega,s)|^2 \leq \frac{C
 \beta^2}{1-\beta^2}\frac{|\tilde{g}'_0|^2}{\eta^2},\quad {\rm Re}\, s = \eta > 0,
\end{equation}
for $(\omega, s)$ in the vicinity of the generalized eigenvalue $s_0=\pm
i\,\sqrt{1-\beta^2}\,\omega_0$. In general, the solution becomes unbounded as $|\beta|\to 1$. The
truncation error terms which perturb the numerical solution are therefore amplified by a factor
$1/\sqrt{1-\beta^2}$, which explains the poor accuracy in the numerical solution as $|\beta|\to 1$.

For boundary data $g_0(y,t)$ which have a Laplace-Fourier transform that can be
written as
\[
\tilde{g}_0(\omega,s) = s \tilde{G}(\omega,s),
\]
estimate (\ref{eq:gen-eig-est}) becomes
\[
|\tilde{u}'(0,\omega,s)|^2 \leq \frac{C
  \beta^2}{1-\beta^2}\frac{|s|^2|\tilde{G}|^2}{\eta^2} \approx 
\frac{C \beta^2}{|s_0|^2/\omega_0^2}\frac{|s|^2|\tilde{G}|^2}{\eta^2} = 
C \beta^2 \omega_0^2\frac{|\tilde{G}|^2}{\eta^2},\quad s\to s_0.
\]
Hence the $|\beta|\to 1$ singularity cancels out and the solution is bounded independently
of $\beta$. The factor `$s$' on the Laplace transform side corresponds to a
time-derivative on the un-transformed side. Hence, the solution is bounded independently
of $\beta$ if the boundary forcing can be written as a time-derivative of a function with
bounded Laplace-Fourier transform,
\[
g_0(y,t) = G_{t}(y,t),\quad G(y,0)=0,\quad \left|\int_{y=0}^1\int_{t=0}^\infty e^{-2\pi i\omega y} e^{-st}G(y,t)\,dtdy\right| <
\infty,\ \ \re s \geq 0.
\]
The latter condition is satisfied if $G(y,t)$ is in $L^1$, i.e.,
\begin{equation}\label{eq:L1}
\int_{y=0}^1 \int_{t=0}^\infty |G(y,t)|\,dtdy < \infty.
\end{equation}

To test this theory numerically, we use a homogeneous interior forcing ($F=0$) and homogeneous initial
conditions ($f_1=f_2=0$), homogeneous forcing on the $x=1$ boundary ($g_1=0$), and
consider three different forcing functions on the $x=0$ boundary: $g_0^{(1)}(y,t)=G(y,t)$,
$g_0^{(2)}(y,t)=G_t(y,t)$, and $g_0^{(3)}(y,t) = G_{tt}(y,t)$. Here we choose $G(y,t)$ to
trigger a surface wave:
\[
G(y,t) = u_s(0,y,t)\, e^{-(t/t_0 - 7)^2},\quad t_0 = 0.2,
\]
where $u_s$ is defined by (\ref{eq:surface-wave}). The Gaussian pulse $\exp(-(t/t_0 -
7)^2)$ decays exponentially fast away from its center at $t=7t_0$. For example, it equals
$1.23\times 10^{-4}$ at $t=7t_0\pm 3t_0$, and $5.24\times 10^{-22}$ at $t=7t_0 \pm
7t_0$. The function $G(y,t)$ satisfies (\ref{eq:L1}), so our theory predicts that boundary
forcings $g_0^{(2)}$ and $g_0^{(3)}$ should give solutions that are bounded independently
of $\beta$. However, the time-integral of a Gaussian pulse is the error-function
($\mbox{erf}$), so the boundary forcing $g_0^{(1)}$ does not satisfy (\ref{eq:L1}).

\begin{figure}[htp]
  \begin{center}
    \includegraphics[width=0.9\textwidth]{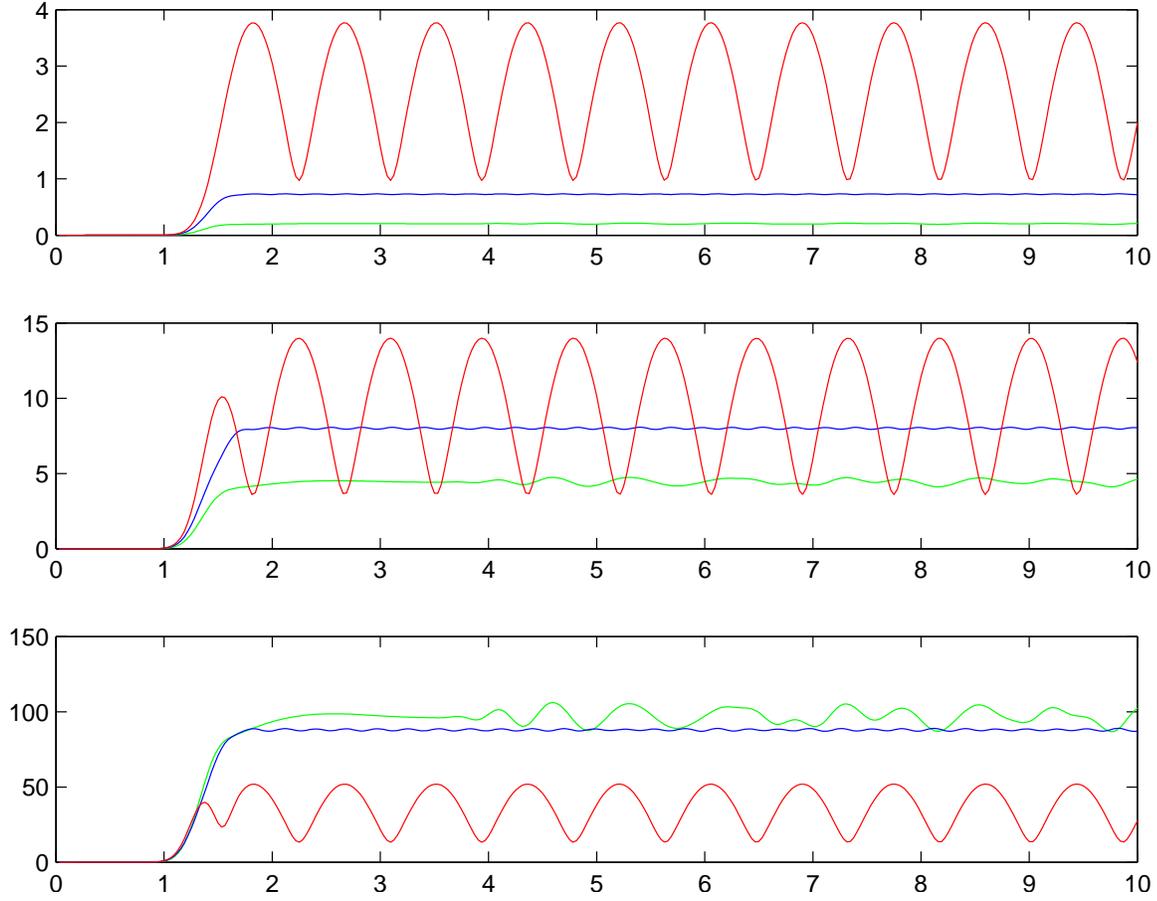}
  \end{center}
  \caption{The max norm of the solution as function of time for the boundary forcing
  functions $g_0^{(1)} = G$ (top), $g_0^{(2)}=G_t$ (middle), and $g_0^{(3)} = G_{tt}$
  (bottom). In each figure, the green, blue, and red curves correspond to $\beta=0.5$,
  $\beta=0.9$, and $\beta=0.99$, respectively.}
  \label{fig:derivatives}
\end{figure}

\begin{figure}[htp]
  \begin{center}
    \includegraphics[width=0.9\textwidth]{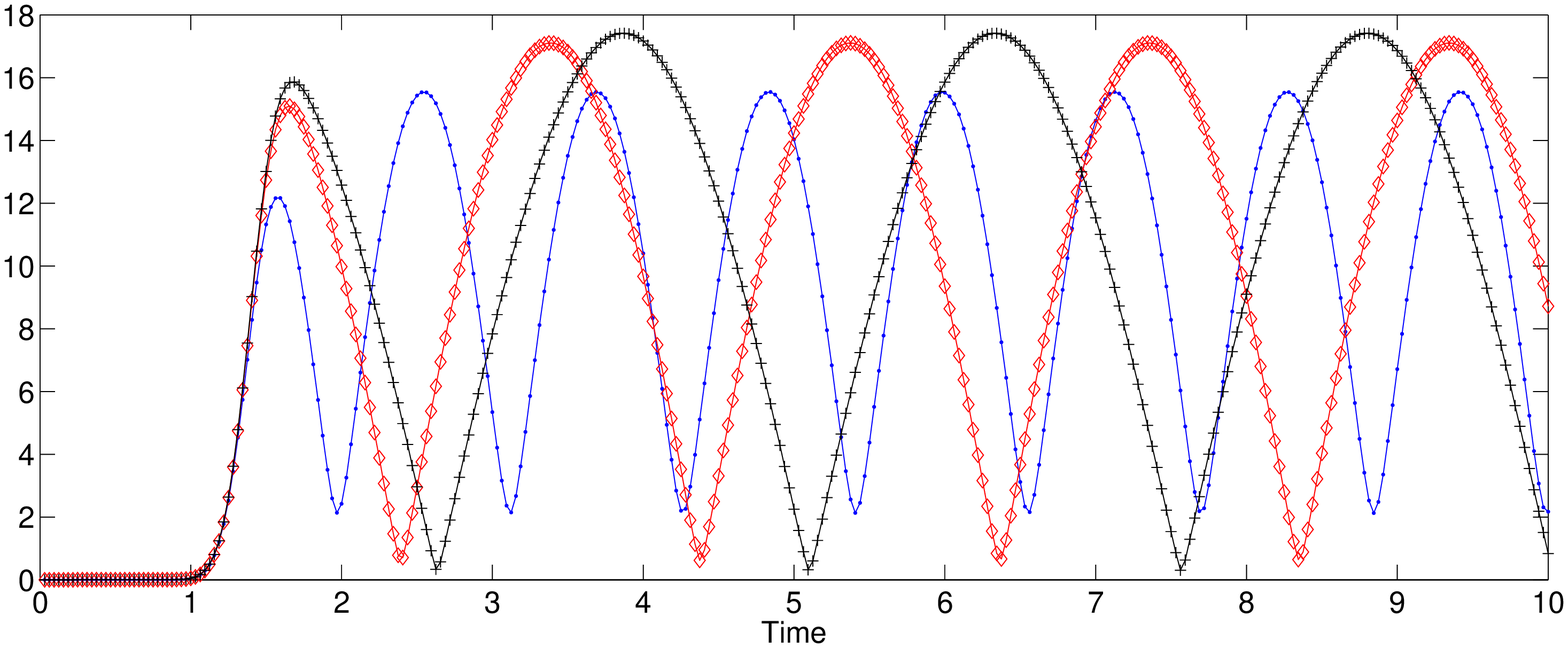}
  \end{center}
  \caption{The max norm of the solution as function of time for the boundary forcing
  function $g_0^{(2)}=G_t$ for $\beta=0.995$ (blue/dots), $\beta=0.999$ (red/diamonds) and
  $\beta=0.9997$ (black/plusses).}
  \label{fig:one-derivative}
\end{figure}

In the numerical calculations we take $\omega_0=8\pi$ and study the cases $\beta=0.5$,
$\beta=0.9$, and $\beta=0.99$. The grid size and time step are $h=2.5\times 10^{-3}$
and $\delta_t=0.5 h$. The max norm of the solution as function of time is shown
in Figure~\ref{fig:derivatives}.
The case $g_0^{(1)}=G$ in the top sub-figure illustrates the general case where the
solution grows as $|\beta|\to 1$. Note that estimate (\ref{eq:gen-eig-est}) predicts the
solution to grow like $\beta/\sqrt{1-\beta^2}$, which means that the solution should be
about 3.5 times larger for $\beta=0.99$ than $\beta=0.9$. In the numerical calculation,
the max norm of the solution grows from about 0.75 for $\beta=0.9$ to 3.75 for
$\beta=0.99$, which is slightly faster than predicted by theory. The case
$g_0^{(3)}=G_{tt}$ in the bottom sub-figure shows the opposite situation when the solution
decays as $\beta\to 1$ because the forcing function is a second time-derivative of a
function with bounded $L^1$ norm, corresponding to an $s^2$ factor on the Laplace transform side. The
intermediate case $g_0^{(2)} = G_t$ is shown in the middle sub-figure. Here the solution
grows between $\beta=0.9$ and $\beta=0.99$, but not as fast as for $g_0^{(1)}$. To more
closely study the behavior near $\beta=1$, we take $\beta=0.995$, 0.999 and 0.9997
corresponding to $\sqrt{1-\beta^2} \approx 0.0998$, 0.0447 and 0.0244, respectively. To
properly resolve the solution we here use an extra fine grid with $h=1.25\times 10^{-3}$ and
$\delta_t = 0.5 h$. The max norm of the solutions, shown in
Figure~\ref{fig:one-derivative}, reveal that the solution indeed stays bounded
independently of $\beta$, confirming our theory also for boundary forcing $g_0^{(2)}=G_t$.

%%%%%%%%%%%%%%%%%%%%%%%%%%%%%%%%%%%%%%%%%%%%%%%%%%%%%%%%%%%%%%%%%%%%%%%%
% old appendix.tex
%%%%%%%%%%%%%%%%%%%%%%%%%%%%%%%%%%%%%%%%%%%%%%%%%%%%%%%%%%%%%%%%%%%%%%%%
\setcounter{thrm}{0}
\renewcommand{\thethrm}{A.\arabic{thrm}}
\section*{Appendix}
In this appendix we collect a number of auxilary lemmas.
\begin{lmm}\label{lemmaA.1}
The solution of
\begin{equation}
y_x=\lambda y+F,\quad {\rm Re}~\lambda >0,\quad 0\le x<\infty \label{A.1}
\end{equation}
satisfies the estimate
\[
|y(0)|^2\le \frac{1}{2 {\rm Re}~\lambda}~\|F\|^2,\quad \|y\|^2\le
\frac{1}{({\rm Re}~\lambda)^2}~\|F\|^2,\quad \|F\|^2=\int_0^\infty |F|^2 dx.
\]
\end{lmm}
\begin{proof} Integration by parts gives us
\[
(y,y_x)=-|y(0)|^2-(y_x,y),\quad {\rm i.e.,}\quad 2{\rm Re}(y,y_x)=-|y(0)|^2.
\]
Therefore
\[\begin{split}
\frac{1}{2}|y(0)|^2+{\rm Re}~\lambda~\|y\|^2=-{\rm Re}~(y,F)&\le \|y\|\,\|F\|\cr
		&\le \frac{\alpha}{2}~{\rm Re}~\lambda~\|y\|^2 +
\frac{1}{2\alpha}\frac{\|F\|^2}{{\rm Re}~\lambda}, \quad \alpha>0.
\end{split}
\]
With $\alpha=2$ the first inequality follows. With $\alpha=1$ the second
inequality follows.
\end{proof}

\begin{lmm}\label{lemmaA.2}
The solution of
\[
y_x=-\lambda y+F,\quad y(0)=g,\quad {\rm Re}~\lambda>0,\quad 0\le x <\infty,
\]
satisfies
\begin{equation}
\|y\|^2\le \frac{1}{{\rm Re}\,\lambda}~|g|^2 + \frac{1}{({\rm
Re}\,\lambda)^2}~\|F\|^2. \label{A.2}
\end{equation}
\end{lmm}
\begin{proof}
\[
\langle y,y\rangle_x=2{\rm Re}\,\langle y,y_x\rangle =-2({\rm Re}\,\lambda)
|y|^2 + 2{\rm Re}\,\langle y,F\rangle.
\]
As $y\in \xLtwo,$ integrating we have
\[\begin{split}
-|y(0)|^2&=-2{\rm Re}\,\lambda~\|y\|^2 + 2 {\rm Re}\,(y,F)\cr
&\le -2{\rm Re}\,\lambda~\|y\|^2 + 2~\|y\|~\|F\|\cr
&\le -2{\rm Re}\,\lambda~\|y\|^2 + {\rm Re}\,\lambda~\|y\|^2 +
\frac{\|F\|^2}{{\rm Re}\,\lambda}.\cr
\end{split}
\]
Thus,
\[
{\rm Re}\,\lambda~\|y\|^2 \le |y(0)|^2 + \frac{\|F\|^2}{{\rm Re}\,\lambda}
\]
and the lemma follows.
\end{proof}

\begin{lmm}\label{lemmaA.3}
Let $a,b$ be real and consider
$ \sqrt{a+ib} $ with $-\pi< {\rm arg}(a+ib)\le \pi,~{\rm arg}\sqrt{z}=\frac{1}{2}
{\rm arg}\, z.$ Then, the following inequalities hold
\begin{equation}
2^{-1/4}\sqrt{|a|+|b|} \le|\sqrt{a+ib}|\le \sqrt{|a|+|b|} \label{A.3}
\end{equation}
\begin{equation}
2^{-3/4}|\sqrt{a+ib}|\le 2^{-3/4}\sqrt{|a|+|b|}\le {\rm Re}\sqrt{a+ib}\le
|\sqrt{a+ib}|\quad {\rm for} \quad a\ge 0, \label{A.4}
\end{equation}
\begin{equation}
2^{-5/4}\frac{|b|}{|\sqrt{a+ib}|}\le 2^{-1} \frac{|b|}{\sqrt{|a|+|b|}}\le {\rm
Re}\sqrt{a+ib}\le |\sqrt{a+ib}|\quad {\rm for} \quad a\le 0. \label{A.5}
\end{equation}
\end{lmm}
\begin{proof} In polar notation $a+ib = \rho e^{i\theta},$ $\rho=\sqrt{a^2+b^2}>0,$
$-\pi<\theta\le\pi,$ and 
\[
\sqrt{a+ib}=\sqrt{\rho}\, e^{i\frac{\theta}{2}}
\]
We have
\[
\sqrt{|a|+|b|}=\sqrt{\rho}\sqrt{|\cos\theta|+|\sin\theta|}\ge
2^{1/4}\sqrt{\rho}
\]
and the first inequality in (\ref{A.3}) follows. The second inequality in
(\ref{A.3}) follows
from the triangle inequality. If $a\ge0$ then $\frac{\theta}{2}\in
[-\pi/4,\pi/4]$ and (\ref{A.3}) implies
\[
|\sqrt{a+ib}|\ge \re \sqrt{a+ib}=\sqrt{\rho}\cos(\theta/2)\ge
\frac{\sqrt{2}}{2}\rho \ge 2^{-3/4}\sqrt{|a|+|b|}\ge 2^{-3/4} |\sqrt{a+ib}|
\]
which is (\ref{A.4}). To prove (\ref{A.5}) notice that, as
$a\le 0,$
\[
\frac{|b|}{\sqrt{|a|+|b|}}\le \frac{|b|}{|\sqrt{a+ib}|}=
\frac{\rho\sin\theta}{\sqrt{\rho}}=\sqrt{\rho}\,2\,
\sin(\theta/2)\cos(\theta/2)\le2\,\re\sqrt{a+ib}
\]
therefore
\[
\frac{1}{2}\frac{|b|}{\sqrt{|a|+|b|}}\le\re\sqrt{a+ib}\le|\sqrt{a+ib}|
\]
and the inequality follows from (\ref{A.3}).
\end{proof}

For the forthcoming lemmata, we remind the reader that $s=\eta+i\xi$, where $\eta$, $\xi\in {\mathbb
  R}$, and $\kappa=\sqrt{\omega^2 + s^2}$. We shall now apply the last lemma to
\[
\kappa=\sqrt{\omega^2+\eta^2-\xi^2+2i\xi\eta},\quad {\rm i.e.,} \quad
a=\omega^2+\eta^2-\xi^2,~b=2\xi\eta.
\]
In the following three lemmas we denote by $\delta$ a constant with $0<\delta <1$.

\begin{lmm}\label{lemmaA.4}
Let
\[
\delta_1=2^{-1/4}\sqrt{\delta},\quad \delta_2=2^{1/2}(1-\delta)^{1/4},\quad
\delta_3= \min(\delta_1,\delta_2). 
\]
Then
\begin{equation}\label{A.6}
|\kappa|\ge \begin{cases}
\delta_1\sqrt{\omega^2+|s|^2} & \mbox{if } |\omega^2 + \eta^2-\xi^2| \ge \delta (\omega^2
+|s|^2)\cr
\delta_2\sqrt{\sqrt{|s|^2+\omega^2} \,\eta} & \mbox{otherwise.}\cr\end{cases}
\end{equation}
Also, always
\begin{equation}
|\kappa|\ge \delta_3\eta. \label{A.7}
\end{equation}
\end{lmm}
\begin{proof} By (\ref{A.3}) we obtain, for the first case,
\[
|\kappa|\ge 2^{-1/4}\sqrt{|\omega^2+\eta^2-\xi^2| +2|\xi|\eta} \ge
2^{-1/4}\sqrt{\delta}\sqrt{\omega^2+|s|^2}.
\]
If $|\omega^2+\eta^2-\xi^2|<\delta(\omega^2+|s|^2),$ then
\[
2\xi^2\ge (1-\delta)(\omega^2+\xi^2+\eta^2)= (1-\delta)(\omega^2+|s|^2).
\]
Therefore
\[
|\kappa|\ge 2^{-1/4}\sqrt{2|\xi|\eta} \ge \sqrt{2\sqrt{1-\delta}\sqrt{\omega^2+|s|^2}\,\eta}.
\]
Also $\sqrt{\omega^2+|s|^2}\ge\eta$ implies (\ref{A.7}). This proves the lemma.
\end{proof}

\begin{lmm}\label{lemmaA.5}
\begin{equation}\label{A.8}
\begin{split}
{\rm Re}\,\kappa &\ge \begin{cases}
2^{-5/4}|\kappa| & \mbox{if }~\omega^2+\eta^2 -\xi^2 \ge 0,\nonumber\\
2^{-3/4} \frac{\sqrt{\omega^2+|s|^2}}{|\kappa|} \eta & \mbox{if }~\omega^2+\eta^2
-\xi^2 < 0.\cr\end{cases}\cr
{\rm Re}\,\kappa &\ge \delta_4\eta, \quad \delta_4=2^{-3/4}\min(1,\delta_3).
\end{split}
\end{equation}
\end{lmm}
\begin{proof}
If $\omega^2+\eta^2-\xi^2\ge 0,$ then (\ref{A.4}) gives us
\[
{\rm Re}\,\kappa \ge 2^{-3/4}|\kappa|\ge 2^{-3/4}\delta_3 \eta.
\]
 If $\omega^2+\eta^2-\xi^2< 0,$ then $2\xi^2\ge 
 \omega^2+\eta^2+\xi^2=\omega^2+|s|^2.$ Therefore (\ref{A.5}) gives us
\[
{\rm Re}\,\kappa \ge 2^{-5/4} \frac{2|\xi|\eta}{|\kappa|}\ge
2^{-3/4} \frac{\sqrt{\omega^2+|s|^2}}{|\kappa|}\,\eta \ge 2^{-3/4}\eta.
\]
This proves the lemma.
\end{proof}

\begin{lmm}\label{lemmaA.6}
\[
|\kappa|\,{\rm Re}\,\kappa \ge \delta_6\sqrt{\omega^2+|s|^2}\,\eta,\quad
\delta_6=\min\Bigl(\delta_1 \delta_4, 2^{-5/4}\delta_2^2, 2^{-3/4}\Bigr).
\]
\end{lmm}
\begin{proof} If $\omega^2+\eta^2-\xi^2\ge 0,$ and
$|\omega^2+\eta^2-\xi^2|\ge \delta(\omega^2+|s|^2),$ then (\ref{A.6}) and
(\ref{A.8}) give us
\[
|\kappa|\,{\rm Re}\,\kappa \ge \delta_1 \sqrt{\omega^2 + |s|^2}~\delta_4 \eta.
\]
If $\omega^2+\eta^2-\xi^2\ge 0,$ and $|\omega^2+\eta^2-\xi^2|< \delta(\omega^2+|s|^2),$ 
then (\ref{A.6}) and (\ref{A.8}) give us
\[
|\kappa|\,\re\kappa \ge 2^{-5/4} |\kappa|^2\ge 2^{-5/4}\delta_2^2 \sqrt{\omega^2+|s|^2}\,\eta.
\]
Finally, if $\omega^2+\eta^2-\xi^2< 0,$ by (\ref{A.8}) we obtain
\[
|\kappa|\,\re\kappa \ge 2^{-3/4} \sqrt{|\omega|^2+|s|^2}\,\eta.
\]
This proves the lemma.
\end{proof}

\begin{lmm}\label{lemmaA.7}
Assume that, for the boundary condition 1), $a>0,~|b|<1.$
Then there is a constant $\delta >0$ such that, for all $\omega$ and
$s$ with ${\rm Re}\,s\ge 0,$
\[
|s+a\kappa -ib\omega|\ge \delta \sqrt{|s|^2+|\omega|^2}.
\]
\end{lmm}
For the proof, see Lemma 3 of \cite{Kreiss2006}.

Finally we have a lemma similar to Lemma A.4 and Lemma A.6 but for the
normalized variables
\[
\kappa'=\sqrt{\eta'^2 - \xi'^2 + \omega'^2 + 2i\xi'\omega'}, \quad
\xi'^2+\omega'^2=1,~|\eta'|\ll 1.
\]
By Lemma A.3,
\begin{equation}\label{A.9}
\begin{split}
|\kappa'|&\ge 2^{-1/4}\sqrt{|-\xi'^2+\omega'^2+\eta'^2|+2|\xi'|\,|\eta'|}\cr
{\rm Re}\,\kappa' &\ge 2^{-3/4} \sqrt{|-\xi'^2+\omega'^2 +
\eta'^2|+2|\xi'|\,|\eta'|} \quad {\rm if}~\xi'^2\le \omega'^2+\eta'^2,\cr
{\rm Re}\,\kappa' &\ge \frac{|\xi'|
\eta'}{\sqrt{|-\xi'^2+\omega'^2+\eta'^2|+2|\xi'|\,|\eta'|}} \quad {\rm
if}~\xi'^2
>\omega'^2+\eta'^2.\cr
\end{split}
\end{equation}
\begin{lmm}\label{lemmaA.8}
There is a constants $\delta>0$ such that
\[
|{\rm Re}\,\kappa'|\ge \delta \eta',\quad |\kappa'|\ge 2^{-1/4}\,\eta',
\quad |\kappa'|\;|{\rm Re}\,\kappa'|\ge \delta_6\,\eta'.
\]
where $\delta_6$ is that of Lemma A.6.
\end{lmm}
\begin{proof} From Lemma Lemma A.6, dividing by $\omega^2+\xi^2,$
\[
|\kappa'|~|{\rm Re}\,\kappa'|~=~\frac{|\kappa|~{\rm Re}\,\kappa}{ \omega^2+\xi^2}~\ge~\delta_6
\frac{\sqrt{\omega^2+|s|^2}}{ \sqrt{\omega^2+\xi^2}}
\frac{\eta}{ \sqrt{\omega^2+\xi^2}}~\ge~ \delta_6~\eta'.
\]
Now if $\omega'^2 - \xi'^2 \ge 0,$ then $\omega'^2 - \xi'^2 +\eta'^2\ge \eta'^2,$
and 
\[
\sqrt{|\omega'^2 - \xi'^2 +\eta'^2|+2|\xi'|\eta'}\ge \sqrt{\eta'^2+2|\xi'|\eta'}\ge \eta',
\]
so that, by Lemma A.3,
\[
|\kappa'|\ge 2^{-1/4} \eta', \quad {\rm and}\quad \re\kappa'\ge 2^{-3/4}\eta'.
\]
If $\omega'^2-\xi'^2\le 0$ then $|\xi'^2|\ge 1/2.$ From (\ref{A.9})
\[
|\kappa'|\ge 2^{-1/4} \sqrt{2\sqrt{2}\,\eta'}\ge 2^{1/2}~\eta'.
\]
This last estimate also holds for $|{\rm Re}\,\kappa'|$ when $\omega'^2 - \xi'^2
+\eta'^2\ge 0.$  When $\omega'^2 - \xi'^2 +\eta'^2\le 0$ we just use $|\omega'^2
- \xi'^2 +\eta'^2|+2|\xi'|\,\eta' \le {\rm const.}$ and from (\ref{A.9}) the estimate
for $\re\kappa'$ follows from (\ref{A.9}). This proves the lemma.
\end{proof}

\begin{acknowledgement}
We thank Prof. Oscar Reula who participated in the original discussions leading to this
work.
\end{acknowledgement}

%%%%%%%%%%%%%%%%%%%%%%%%%%%%%%%%%%%%%%%%%%%%%%%%%%%%%%%%%%%%%%%%%%%%%%
%       BIBLIOGRAPHY
%%%%%%%%%%%%%%%%%%%%%%%%%%%%%%%%%%%%%%%%%%%%%%%%%%%%%%%%%%%%%%%%%%%%%%


\begin{thebibliography}{1}

\bibitem{Kreiss1970}
H.-O. Kreiss.
\newblock Initial boundary value problems for hyperbolic systems.
\newblock {\em Commun. Pur. Appl. Math.}, 23:277--298, 1970.

\bibitem{Kreiss-Lorenz1989}
H.-O. Kreiss and J.~Lorenz.
\newblock {\em Initial-boundary value problems and the {N}avier-{S}tokes
  equations}.
\newblock Academic Press, San Diego, 1989.

\bibitem{Gustafsson-Kreiss-Oliger1995}
B. Gustafsson, H.-O. Kreiss and J. Oliger.
\newblock {\em Time dependent problems and difference methods}.
\newblock Wiley-Interscience, 1995.

\bibitem{Kreiss2006}
H.-O. Kreiss and J.~Winicour.
\newblock Problems which are well posed in a generalized sense with
  applications to the {E}instein equations.
\newblock {\em Class. Quantum Grav.}, 23(16):405--420, 2006.

\bibitem{Agranovich1972}
M.~S. Agranovich.
\newblock Theorem on matrices depending on parameters and its applications to
  hyperbolic systems.
\newblock {\em Funct. Anal. Appl.}, 6:85--93, 1972.

\end{thebibliography}
\end{document}